\documentclass[]{article}

\usepackage{amsfonts,amsmath,amssymb}
\usepackage{amsthm}
\usepackage[cp866]{inputenc}
\usepackage[english]{babel}
\usepackage[hyperindex,colorlinks,citecolor=blue,pagebackref]{hyperref}

\allowdisplaybreaks

\theoremstyle{plain}
\newtheorem{theorem}{Theorem}[section]
\newtheorem{corollary}[theorem]{Corollary}
\newtheorem{proposition}[theorem]{Proposition}

\newtheorem{lemma}[theorem]{Lemma}

\newtheorem{fact}[theorem]{Fact}
\newtheorem{question}{Question}

\theoremstyle{definition}
\newtheorem{definition}[theorem]{Definition}

\theoremstyle{remark}
\newtheorem{claim}{Claim}
\newtheorem{remark}[theorem]{Remark}

\newcommand{\ZZ}{\mathbb{Z}}

\newcommand{\I}{\mathbf{I}}
\newcommand{\R}{\mathbb{R}}
\newcommand{\N}{\mathbb{N}}

\renewcommand{\d}{\mathbf{d}}
\newcommand{\e}{\mathbf{e}}
\newcommand{\hil}{Q}
\newcommand{\T}{\mathbb{T}}
\newcommand{\U}{\mathcal{U}}
\newcommand{\V}{\mathcal{V}}
\newcommand{\A}{\mathcal{A}}
\newcommand{\B}{\mathcal{B}}
\newcommand{\F}{\mathcal{F}}
\newcommand{\X}{\mathcal{X}}
\newcommand{\x}{\mathbf{x}}
\newcommand{\zero}{\mathbf{0}}
\newcommand{\K}{\mathcal{K}}

\renewcommand{\P}{\mathcal{P}}
\newcommand{\Baire}{\N^\N}

\newcommand{\Sph}{\mathbf{S}}
\newcommand{\im}{\mathop{\mathrm{im}}}
\newcommand{\supp}{\mathop{\mathrm{supp}}}
\newcommand{\cB}{\overline{B}}
\newcommand{\Q}{\mathbb{Q}}
\newcommand{\J}{\mathcal{J}}

\newcommand{\Spec}{\mathrm{Sp}}
\newcommand{\inj}{\mathcal{I}}
\newcommand{\CP}{\mathbf{CP}}

\newcommand{\zd}{\ZZ_2}

\newcommand{\om}{\omega}

\newcommand{\res}[1]{\mathop{\upharpoonright_{#1}}}

\title{Degree Spectra of Homeomorphism Type of Compact Polish Spaces}
\author{Mathieu Hoyrup\thanks{This work benefited of Inria's program Chercheur Invit\'e}\\
Universit\'e de Lorraine, CNRS, Inria, LORIA, France\\
{\tt mathieu.hoyrup@inria.fr}\\
and\\
Takayuki Kihara\thanks{Kihara's research was partially supported by JSPS KAKENHI Grant 19K03602, 15H03634, and the JSPS Core-to-Core Program (A. Advanced Research Networks).}\\
Graduate School of Informatics, Nagoya University, Japan\\
{\tt kihara@i.nagoya-u.ac.jp}\\
and\\
Victor Selivanov\thanks{Selivanov's contributions to most of the results of this paper, in particular to those in Section 3, were supported by the Russian Science Foundation, project 19-71-30002}
\\Department of Mathematics and Computer Science,\\
St. Petersburg State University, Saint Petersburg, Russia, and\\
A.P. Ershov Institute of Informatics Systems SB RAS, Novosibirsk, Russia\\
{\tt vseliv@iis.nsk.su}
}

\begin{document}
\large
\date{}
 \maketitle

\begin{abstract}
A Polish space is not always homeomorphic to a computably presented Polish space.
In this article, we examine degrees of non-computability of presenting homeomorphic copies of compact Polish spaces.
We show that there exists a $\zero'$-computable low$_3$ compact Polish space which is not homeomorphic to a computable one, and that, for any natural number $n\geq 2$, there exists a Polish space $X_n$ such that exactly the high$_{n}$-degrees are required to present the homeomorphism type of $X_n$. Along the way we investigate the computable aspects of \v{C}ech homology groups. We also show that no compact Polish space has a least presentation with respect to Turing reducibility.

The first version of this article appeared in April 2020. A major update was made in September 2023, with improved proofs and results. This is the final version from January 2024, with more results on \v{C}ech homology groups.


{\bf Keywords}: computable topology, computable presentation, computable Polish space, degree spectrum, \v{C}ech homology groups
\end{abstract}

\tableofcontents

\section{Introduction}\label{sec_intro}

How difficult is it to describe an explicit presentation of an abstract mathematical structure?
For a particular structure~$S$, what is the degree spectrum of its isomorphism type, i.e.~what are the Turing degrees that compute a presentation of a structure that is isomorphic to~$S$? For a class of structures, what are the possible degree spectra of their isomorphism types? These have long been the fundamental questions in computable structure theory, and researchers in this area have obtained a huge number of interesting results on Turing degrees of presentations of isomorphism types of groups, rings, fields, linear orders, lattices, Boolean algebras, and so on~\cite{AK00,HRM98,FHM14,HKSS02}.

In this article we focus on presentations of compact Polish spaces. The notion of a presentation plays a central role, not only in computable structure theory, but also in computable analysis \cite{HCCA,BHW08,Wei00}. In this area, one of the most crucial problems was how to present large mathematical objects (which possibly have the cardinality of the continuum) such as metric spaces, topological spaces and so on, and then researchers have obtained a number of answers to this question.
In particular, the notion of a computable presentation of a Polish space can be traced back to the 1960s; e.g.~\cite{Mos64}. Computable Polish spaces have been extensively studied in computable analysis \cite{HCCA,PoRi89,Wei00} and descriptive set theory \cite{Mos09}.

In recent years, several researchers have succeeded in obtaining various results on Turing degrees of {\em isometric isomorphism types} of Polish spaces, separable Banach spaces, and Polish groups; see~\cite{CMS19,McSt19,Mel13,LMN23}. All these works are devoted to metric structures. Results on Turing degrees of \emph{homeomorphism types} of Polish spaces have been obtained only very recently, and independently of the present article, in \cite{HTMN20,Melnikov21,LMN23,DM23}. Computable Polish groups were investigated in \cite{GMNT18,M18}, and general topological spaces were studied in \cite{Se19} in analogy with the earlier investigation of degrees of isomorphism types of algebraic structures. Some results were also obtained for domains.

In this article, a Polish presentation of a Polish space~$X$ is a dense sequence~$(s_i)_{i\in\N}$ in~$X$ and a complete metric~$d$ inducing the topology; for a Turing degree~$\d$, such a presentation is~$\d$-computable if the real numbers~$d(s_i,s_j)$ are uniformly~$\d$-computable. The Polish degree spectrum of~$X$ is then the set of Turing degrees that compute a Polish presentation of~$X$. We will also be interested in compact presentations of a compact Polish space~$X$, that carry a Polish presentation together with an enumeration of the finite open covers of~$X$. The compact degree spectrum of~$X$ is the set of Turing degrees that compute a compact presentation of~$X$.

One of the first questions to be asked about presentations of Polish spaces is:
\begin{question}\label{question:basic}
Does there exist a Polish space with a~$\zero'$-computable Polish presentation but not computable one?
\end{question}
Observe that there are continuum many homeomorphism types of Polish spaces, but there exist only countably many computable presentations of Polish spaces, so there is a Polish space which is not homeomorphic to any computably presented Polish space.

Countable spaces cannot be used to solve this problem because of the ``hyperarithmetic-is-recursive'' phenomenon, see~\cite{GrMo08}; see also Section \ref{sec:Cantor-Bendixson}.
In this article we answer Question \ref{question:basic} in the affirmative. We note that Harrison-Trainor, Melnikov, and Ng \cite{HTMN20} have recently obtained an independent solution, and via a different proof.
One possible approach to solve this problem is using Stone duality between countable Boolean algebras and zero-dimensional compact Polish spaces; see also Section \ref{boolean}.
Combining this idea with classical results on isomorphism types of Boolean algebras \cite{KnSt00}, one can conclude that every {\em low$_4$-presented} zero-dimensional compact Polish space is homeomorphic to a computable one.
This was also independently noticed by Harrison-Trainor, Melnikov, and Ng \cite{HTMN20}.

Our next step is to develop new techniques beyond Stone duality. More specifically, the next question is whether there exists a Polish space whose Polish degree spectrum is different from that of a zero-dimensional compact space. In particular, it is natural to ask the following:

\begin{question}\label{question:basic2}
Does there exist a Polish space with a low$_4$ Polish presentation but no computable one?
\end{question}

One of the main results of this article is that for every Turing degree~$\d$ and every~$n\geq 1$, there exists a compact Polish space~$\X_{\d,n}$ whose compact degree spectrum is~$\{\x:\d\leq \x^{(n)}\}$ and Polish degree spectrum is~$\{\x:\d\leq \x^{(n+1)}\}$ (Corollaries \ref{cor_one_jump} and \ref{cor_iteration}). This result was independently obtained by Melnikov \cite{Melnikov21}, Lupini, Melnikov and Nies \cite{LMN23} and Downey and Melnikov \cite{DM23} using different arguments, namely by proving that \v{C}ech cohomology groups of a compact space are computable. For~$n=1$ the result presented here is obtained by investigating the computable aspects of \v{C}ech \emph{homology} groups. A part of the analysis of homology groups is a direct adaptation of the analysis of cohomology groups from \cite{LMN23,DM23}, then the arguments diverge: \v{C}ech homology groups are not computable in general, but they are sufficiently effective for our purpose. More precisely, we show that the non-triviality of \v{C}ech homology groups with coefficients in~$\ZZ/2\ZZ$ is~$\Sigma^0_2$; along the way, we also show that it is~$\Sigma^1_1$-complete when taking coefficients in~$\ZZ$.

The result gives a solution to Questions \ref{question:basic} and \ref{question:basic2}: there is a space with a~$\zero'$-computable low$_3$ Polish presentation but no computable Polish presentation, namely the space~$\X_{\d'',1}$ for any~$\d\leq\zero'$ which is low$_3$. 
It also implies the existence of a compact Polish space having a computable Polish presentation but no computable compact presentation (Corollary \ref{cor_polish_nocompact}).

Another consequence is that for~$n\geq 2$ there is a space whose Polish degree spectrum is the set of high$_n$ degrees, namely~$\X_{\d,n-1}$ where~$\d=\zero^{(n+1)}$. This also clarifies substantial differences between zero-dimensional compact Polish spaces and infinite dimensional ones since the class of high$_n$-degrees is never the degree spectrum of a Boolean algebra \cite{JoSo94}.

Another important question is whether a given Polish space has a least Turing degree in its Polish degree spectrum. For instance, it is known that the isomorphism types of linear orders, trees, abelian $p$-groups, etc.~have no least presentation whenever they are not computably presentable, see~\cite{FHM14}.

\begin{question}\label{question:basic3}
Does there exist a compact Polish space with no computable Polish presentation, but whose Polish degree spectrum contains a least Turing degree?
\end{question}

We answer Question \ref{question:basic3} in the negative.
More precisely, we show the cone-avoidance theorem for compact Polish spaces, which states that, for any non-c.e.~set $A\subseteq\N$, every compact Polish space has a Polish presentation that does not enumerate $A$ (Theorem \ref{thm_cone}).

It contrasts with a result proved in \cite{LMN23} that for each set~$A\subseteq\N$, there exists a compact Polish space whose \emph{compact} degree spectrum is~$\{\x:A\text{ is $\x$-c.e.}\}$. If~$(p_i)_{i\in\N}$ is the increasing enumeration of the prime numbers, this space is defined as the Pontryagin dual~$\hat{G}$ of the subgroup~$G$ of~$\Q$ generated by the elements~$\frac{1}{p_i}$ with~$i\in A$.

We finally prove a precise relationship between the Polish degree spectrum and the compact degree spectrum of compact Polish spaces that are \emph{perfect} (Corollary \ref{cor_degree_perfect}), and show that it fails for some non-perfect space (Proposition \ref{prop_nonperfect}).

The article is organized as follows. In Section \ref{sec_prelim} we present the necessary background and develop the technical tools that will be used throughout the article. In Section \ref{sec:Cantor-Bendixson} we briefly discuss the effective aspects of the Cantor-Bendixson derivative and of Stone duality. In Section \ref{sec_homo} we investigate the computable aspects of the \v{C}ech homology groups, which will be used in Section \ref{sec_spectra}, where we show how to realize certain sets of Turing degrees as degree spectra of compact Polish spaces. In Section \ref{sec_cone} we show the cone-avoidance theorem. In Section \ref{sec:jump} we compare degree spectra and compact degree spectra.  

\section{Preliminaries}\label{sec_prelim}

Basic terminology and results of computability theory and computable structure theory are summarized in \cite{AK00}. For basics of computable analysis, we refer the reader to \cite{AvBr14,BHW08,HCCA,Wei00}. For the basic definitions and facts of general topology and dimension theory, see~\cite{HurWal41,vM01}.

\subsection{Presentations of Polish spaces}

If~$X$ is a Polish space, then a {\em Polish presentation} (or simply a {\em presentation}) of~$X$ is~$(X,d,S)$ where~$S=(s_n)_{n\in\N}$ is a dense sequence in~$X$ and~$d$ is a complete metric which is compatible with the topology. For a discussion of presentations of Polish spaces, see also \cite{GKP17}.

For a Turing degree~$\d$, a presentation~$(X,d,S)$ of~$X$ is~$\d$-computable if the real numbers~$d(s_i,s_j)$ are uniformly~$\d$-computable.  A {\em $\mathbf{d}$-computable Polish space} is a Polish space which has a $\mathbf{d}$-computable presentation. For a Polish space $X$, its {\em Polish degree spectrum} $\Spec(X)$ is the set of all Turing degrees $\mathbf{d}$ such that $X$ has a $\mathbf{d}$-computable Polish presentation.

Let~$(X,d,S)$ be a presentation of a Polish space~$X$. A \emph{rational open ball} is~$B(s_i,r)=\{x\in X:d(s_i,x)<r\}$ where~$r>0$ is rational. A \emph{rational open set} is a finite union of rational open balls. A code of a finite rational open cover of $X$ is a finite set~$E\subseteq\N\times\mathbb{Q}_{>0}$ such that for any $X=\bigcup_{(i,r)\in E}B(s_i,r)$. If $X$ is compact, then a {\em compact presentation} of $X$ is a presentation of $X$ equipped with an enumeration of the codes of all finite rational open covers of $X$. In particular, a compact presentation contains an information of total boundedness; that is, a function $\ell:\N\to\N$ such that for all~$n$, $\{B(a_i,2^{-n}):i<\ell(n)\}$ covers the whole space $X$.

A compact Polish space~$X$ is {\em $\mathbf{d}$-computably compact} if it has a $\mathbf{d}$-computable compact presentation. Its {\em compact degree spectrum} $\Spec_c(X)$ is the set of all Turing degrees $\mathbf{d}$ such that $X$ has a $\mathbf{d}$-computable compact presentation.

We will often use the next elementary result.
\begin{lemma}\label{lem:Polish-to-compact}
Let ${X}$ be a compact Polish space. If ${X}$ has a $\mathbf{d}$-computable Polish presentation, then ${X}$ has a $\mathbf{d}'$-computable compact presentation.
\end{lemma}

\begin{proof}
Assume that~$X$ has a~$\d$-computable Polish presentation~$(X,d,S)$. By compactness of $X$, one can observe that~$E$ is a code of a finite rational open cover of $X$ if and only if there exists~$s\in\N$ such that for all $x\in X$ we have $d(x,s_i)\leq r-2^{-s}$ for some $(i,r)\in E$. The latter is equivalent to the existence of $s\in\N$ such that for all $j\in\N$ we have $d(s_j,s_i)\leq r-2^{-s}$ for some $(i,r)\in E$. As $E$ is finite, this is a $\Sigma^0_2$ condition relative to~$\d$, so it is~$\d'$-c.e. In other words, $X$ has a $\mathbf{d}'$-computable compact presentation.
\end{proof}

In Section \ref{sec:jump}, we will present a sharp analysis of this relationship between Polish and compact presentations.


\paragraph{Hyperspaces.}
Every Polish space embeds in the Hilbert cube, and this fact induces an equivalent definition of Polish and compact presentations. The Hilbert cube~$\hil=[0,1]^\N$ is endowed with the complete metric~$d_Q(x,y)=\sum_i2^{-i}|x_i-y_i|$. A point of~$Q$ is rational if its coordinates are rational and only finitely many of them are non-zero. The rational points of~$Q$, enumerated in a canonical way, make~$\hil$ a computable Polish space which is computably compact.

Let~$\V(\hil)$ be the hyperspace of compact subsets of~$\hil$ endowed with the lower Vietoris topology. A subbasis is given by~$\{K\subseteq\hil:K\cap B\neq\emptyset\}$, where~$B$ is a rational ball in~$\hil$ (technically, we need to add~$\{\emptyset\}$, which is the singleton containing the empty compact set, to make it a subbasis).

Let~$\K(\hil)$ be the hyperspace of compact subsets of~$\hil$ endowed with the Vietoris topology. A subbasis is given by a subbasis for the lower Vietoris topology, together with~$\{K\subseteq\hil:K\subseteq U\}$, where~$U$ is a  rational open set in~$\hil$. The Hausdorff metric~$d_H$ is a complete metric generating the Vietoris topology, and the dense sequence of finite sets of rational points of~$\hil$ makes~$\K(\hil)$ a computable Polish space which is computably compact.

We say that a compact set~$K\subseteq Q$ is \emph{$\d$-computably overt} if~$\d$ computes an enumeration of the basic neighborhoods of~$K$ in the lower Vietoris topology. Equivalently,~$K$ is~$\d$-computably overt if it contains a dense computable sequence. We say that~$K$ is \emph{$\d$-computably compact} if~$\d$ computes an enumeration of the basic neighborhoods of~$K$ in the Vietoris topology. Equivalently,~$K$ is~$\d$-computably compact if~$K$ is~$\d$-computably overt and~$K\in\Pi^0_1(\d)$, or if~$\d$ computes a sequence of finite sets~$K_n\subseteq\hil$ of rationals points such that~$d_H(K_n,K)<2^{-n}$.

The next result is folklore. The two parts appear in \cite[Fact 2.11]{AH23} and \cite[Theorem 3.36]{DM23} respectively.

\begin{proposition}\label{prop:presen-represen}
A compact Polish space ${X}$ has a $\mathbf{d}$-computable Polish presentation if and only if it has a~$\d$-computably overt copy in~$Q$.

A compact Polish space ${X}$ has a $\mathbf{d}$-computable compact presentation if and only if it has a~$\d$-computably compact copy in~$\hil$.
\end{proposition}

\begin{proof}
It is essentially an effective version of the fact that every Polish space embeds in~$\hil$.

Let~$(X,d,S)$ be a~$\d$-computable Polish presentation of~$X$. Consider the function~$f:X\to \hil$ mapping~$x\in X$ to~$(x_i)_{i\in\N}\in\hil$ defined by~$x_i=d(x,s_i)/(1+d(x,s_i))$. It is computable and one-to-one, so it is a homeomorphism as~$X$ is compact. Its image is~$\d$-computably overt, because it contains the dense~$\d$-computable sequence~$(f(s_n))_{n\in\N}$.

If~$(X,d,S)$ is moreover~$\d$-computably compact, then so is~$f(X)$, because the image of a~$\d$-computably compact set by a computable function is~$\d$-computably compact.

Conversely, if a compact set~$K\subseteq Q$ is~$\d$-computably overt, then it has a~$\d$-computable Polish presentation, using the metric~$d_Q$ and a dense~$\d$-computable sequence in~$K$. If~$K$ is moreover~$\d$-computably compact, then this presentation is~$\d$-computably compact.
\end{proof}
\subsection{Realizations}
We recall usual operations on spaces, such as disjoint unions and wedge sum, and show that they preserve computability notions.

We will implicitly use the fact that~$[0,1]^n\times Q$ and~$Q\times Q$ are computably homeomorphic to~$Q$, so the results of the next constructions are subsets of~$Q$.
\begin{itemize}
\item If~$X,Y\subseteq Q$ then their disjoint union is~$X\amalg Y=(\{0\}\times X)\cup (\{1\}\times Y)$,
\item If~$X,Y\subseteq Q$ both contain~$\overline{0}$, then their wedge sum is~$X\vee Y=(X\times \{\overline{0}\})\cup (\{\overline{0}\}\times Y)$,
\item If~$X_n\subseteq Q$ for all~$n\in\N$, then their disjoint union is~$\coprod_nX_n=\bigcup_n \{0\}^n\times\{1\}\times X_n$, and its one-point compactification is~$\alpha_0(\coprod_nX_n)=\{\overline{0}\}\cup \coprod_nX_n$, where~$\overline{0}=(0,0,0,\ldots)\in Q$.
\end{itemize}

When~$X,Y$ are Polish spaces, their disjoint union~$X\amalg Y$ is implicitly defined as the Polish space obtained by embedding~$X$ and~$Y$ in~$Q$ and applying the previous definition, and similarly for the wedge sum and the one-point compactification of the disjoint union.
\begin{proposition}\label{prop_constructs}
If~$X,Y,X_n$ are uniformly~$\d$-computably overt (resp.~compact), then~$X\amalg Y,X\vee Y$ and~$\alpha_0(\coprod_n X_n)$ are~$\d$-computably overt (resp.~compact).
\end{proposition}
\begin{proof}
The result is straightforward for finite disjoint unions and wedge sums, because  the sets~$\{0\}\times X$, $\{1\}\times Y$, $X\times \{\overline{0}\}$ and $\{\overline{0}\}\times Y$ easily inherit the computability properties of~$X$ and~$Y$, and so do their unions.

Let us consider the one-point compactification of the disjoint union.  The function~$f_n:Q\to Q$ sending~$x=(x_0,x_1,\ldots)$ to~$(0,\ldots,0,1,x_0,x_1,\ldots)$, starting with~$n$ occurrences of~$0$, is computable uniformly in~$n$. Equivalently, the preimages of rational balls by~$f_n$ are effectively open, uniformly in~$n$.

Assume that~$X_n$ are uniformly~$\d$-computably overt. A rational open ball~$B$ intersects the set~$\alpha_0(\coprod_n X_n)$ iff~$f_n^{-1}(B)$ intersects~$X_n$ for some~$n$, which is~$\d$-c.e.

Now assume that~$X_n$ are uniformly~$\d$-computably compact. A rational open set~$U$ contains~$\alpha_0(\coprod_nX_n)$ iff there exists~$n$ such that~$\{0\}^n\times Q\subseteq U$ and~$X_i\subseteq f_i^{-1}(U)$ for all~$i<n$, which is~$\d$-c.e.
\end{proof}

The next result will be a building block for constructing spaces that encode information about a set of natural numbers.
\begin{proposition}\label{prop_sequence_spaces}
Let~$X_\infty$ and~$(X_n)_{n\in\N}$ be uniformly computably compact subsets of~$Q$ satisfying the following conditions:
\begin{itemize}
\item For all~$n\in\N$, $X_n\subseteq X_{n+1}\subseteq X_\infty$,
\item $d_H(X_n,X_\infty)<2^{-n}$.
\end{itemize}
To a set~$A\subseteq\N$ we associate the compact set
\begin{equation*}
X_A=\begin{cases}
X_{\min A}&\text{if~$A\neq\emptyset$},\\
X_\infty&\text{if~$A=\emptyset$}.
\end{cases}
\end{equation*}

The set~$X_A$ is uniformly computably compact relative to~$A$ and uniformly computably overt relative to any enumeration of~$\N\setminus A$.
\end{proposition}
\begin{proof}
For~$s\in\N$ let~$X_A[s]=X_s$ if~$A\cap [0,s]$ is empty, and~$X_A[s]=X_n$ where~$n=\min(A\cap [0,s])$ otherwise. The sequence~$(X_A[s])_{s\in\N}$ can be computed from~$A$ and satisfies~$d_H(X_A[s],X_A)<2^{-s+1}$, so~$X_A$ is computably compact relative to~$A$.

We show that a rational ball~$B$ intersects~$X_A$ if and only if there exists~$n\in\N$ such that~$[0,n-1]\cap A=\emptyset$ and~$B$ intersects~$X_n$. This condition is c.e.~relative to any enumeration of~$\N\setminus A$. If~$[0,n-1]\cap A=\emptyset$ and~$B$ intersects~$X_n$, then~$X_A$ contains~$X_n$, so~$B$ intersects~$X_A$. Conversely, if~$B$ intersects~$X_A$ then either~$A=\emptyset$ and any sufficiently large~$n$ satisfies the conditions, or~$A\neq\emptyset$ and~$n=\min A$ satisfies the conditions.
\end{proof}

\subsection{Good covers}
We will often extract information about a compact Polish space from its finite open covers. In order for this extraction to be computable, we need to decide which open sets of such a cover intersect, which is made possible by considering good covers only. The content of this section is essentially folklore. Closely related results have appeared in \cite{YasMorTsu99,BraPre03,Kamo05,Iljazovic10,PaulyS020,DM23}.

Let~$(X,d,S)$ be a Polish presentation of~$X$. To a rational open ball~$B(s_i,r)$ we associate the corresponding \emph{rational closed ball}~$\cB(s_i,r)=\{x\in X:d(s_i,x)\leq r\}$, and to a rational open set~$U$ we associate the corresponding union~$\overline{U}$ of rational closed balls. The closed ball always contains the closure of the open ball and is a~$\Pi^0_1$-subset of~$X$, relative to the presentation of~$X$. 

A \emph{good open cover} of~$X$ is a family of rational open sets~$\U=(U_i)_{i\in I}$ where~$I\subseteq \N$ is finite, such that~$X=\bigcup_{i\in I}U_i$ and for every~$J\subseteq I$, if~$\bigcap_{i\in J}U_i=\emptyset$ then~$\bigcap_{i\in J}\overline{U}_i=\emptyset$. An open cover~$(U_i)_{i\in I}$ is a \emph{strong refinement} of an open cover~$(V_j)_{j\in J}$ if for every~$i\in I$, there exists~$j\in J$ such that~$\overline{U}_i\subseteq V_j$.

The next result can be found in \cite[Theorems 1.1 and 3.16]{DM23}.

\begin{lemma}\label{lem_good_covers}
Given a compact presentation of~$X$, one can compute a strong refining sequence of good open covers, i.e.~a sequence~$(\U_s)_{s\in\N}$ of good open covers such that~$\U_{s+1}$ strongly refines~$\U_s$.
\end{lemma}
\begin{proof}
The mesh of a finite open cover~$\U=(U_i)_{i\in I}$ is the maximal diameter of the~$U_i$'s. We first show that there are good open covers of arbitrarily small meshes, and that if~$\U$ is a finite open cover of~$X$ then there exists~$\epsilon>0$ such that every finite cover of mesh~$<\epsilon$ strongly refines~$\U$.

Let~$\epsilon>0$ and~$(B(s,\epsilon))_{s\in I}$ be a rational open cover of~$X$, where~$I$ is finite. By compactness, let~$\delta<\epsilon$ be such that~$(B(s,\delta))_{s\in I}$ is still a cover of~$X$. There are only finitely many values of~$r\in [\delta,\epsilon]$ such that~$(B(s,r))_{r\in I}$ is not a good cover; indeed, these values are~$\min_{x\in X}\max_{s\in J}d(x,s)$, for~$J\subseteq I$; indeed, when the cover is not good there exists~$J\subseteq I$ and a point~$x$ at distance~$\leq r$ from every~$s\in J$, but no point~$y$ at distance~$<r$ from every~$s\in J$, so~$\max_{s\in J}d(x,s)=r$ and~$x$ minimizes this quantity. Therefore, there exists a rational number~$r\in [\delta,\epsilon]$ avoiding this finite set of values, providing a good open cover of mesh~$<2\epsilon$.

Let~$\U$ be a finite open cover of~$X$. By compactness,~$\U$ has a Lebesgue number~$\delta>0$, which means that any set of diameter~$\leq\delta$ is contained in some member of~$U$. If~$\V=(V_i)_{i\in J}$ is a good cover of mesh~$\leq \delta$, then each~$\overline{V_i}$ has diameter~$\leq\delta$, so~$\V$ strongly refines~$\U$.

Now we are given a compact presentation of~$X$ as oracle. Whether a finite rational cover is good is~$\Sigma^0_1$, and whether it strongly refines another rational cover is~$\Sigma^0_1$, and whether its mesh is~$<q$ is~$\Sigma^0_1$. Therefore, we start with searching for some good open cover~$\U_0$ of mesh~$<1$ and inductively look for a good open cover~$\U_{s+1}$ of mesh~$<2^{-s}$ strongly refining~$\U_s$. These objects exist as proved above, and can be effectively found by exhaustive search.
\end{proof}

\subsection{Clopen subsets}

We show that from a compact presentation of~$X$, one can compute a presentation of the boolean algebra of clopen subsets.

\begin{proposition}\label{prop_clopen}
Given a compact presentation of~$X$, one can compute an enumeration~$(C_i)_{i\in\N}$ of the clopen subsets of~$X$, such that equality, inclusion and the finite boolean operations are computable.
\end{proposition}

\begin{proof}
Assume that a compact presentation is given as oracle. One can compute an enumeration~$(U_i,V_i)_{i\in\N}$ of all the pairs of rational open sets satisfying~$X=U_i\cup V_i=X$ and~$\overline{U}_i\cap\overline{V}_i=\emptyset$, because these conditions are c.e.~relative to a compact presentation of~$X$. Let then~$C_i=U_i$ for each~$i$. One can then easily compute the boolean operations: given~$i$, one can compute~$j$ such that~$(U_j,V_j)=(V_i,U_i)$ so~$C_j=X\setminus C_i$; given~$i,j$, one can compute~$k$ such that~$(U_k,V_k)=(U_i\cap U_j,V_i\cup V_j)$, so that~$C_k=C_i\cap C_j$, and symmetrically for the union. The set~$\{i\in\N:C_i=\emptyset\}$ can be computed, because~$C_i=\emptyset\iff X\subseteq V_i\iff U_i=\emptyset$, which is both c.e.~and co-c.e.~relative to a compact presentation. Therefore, equality and inclusion are decidable because they reduce to emptiness of some boolean combination.
\end{proof}

In particular, whether~$X$ is connected is~$\Pi^0_1$ relative to a compact presentation of~$X$.

\subsection{Covering dimension}

Let us briefly recall the covering dimension of a topological space~$X$. If~$\U,\V$ are two open covers of~$X$, we say that~$\V$ is a refinement of~$\U$ if every~$V\in \V$ is contained in some~$U\in\U$. An open cover~$\U$ has \emph{order}~$n$ if the intersection of any~$n$ elements of~$\U$ is empty. The \emph{covering dimension} of~$X$, written~$\dim(X)$, is the least number~$n$ such that every open cover~$\U$ of~$X$ has a refinement~$\V$ of order~$n+1$, if it exists.

\begin{proposition}\label{prop_dim_pi02}
Let~$n\in\N$. Given a compact presentation of~$X$, the predicate~$\dim X\leq n$ is~$\Pi^0_2$, uniformly in~$n$.
\end{proposition}
\begin{proof}
We are given a compact presentation of~$X$ as oracle.
As~$X$ is compact, it is routine to check that in the definition of dimension, one can assume that~$\U$ and~$\V$ are good covers and that~$\V$ strictly refines~$\U$. Therefore,~$\dim X\leq n$ iff for every good cover~$\U$, there exists a good cover~$\V$ that strongly refines~$\U$, and such that for all~$V_0,\ldots,V_n\in\V$, one has~$\overline{V}_0\cap\ldots \cap \overline{V}_n=\emptyset$. This is a~$\Pi^0_2$ predicate, because one can computable enumerate the good covers of~$X$, and the strong refinements of a given good cover.
\end{proof}

\begin{proposition}\label{prop_order}
If~$\dim(X)\leq n$, then given a compact presentation of~$X$, one can compute a strong refining sequence of good open covers~$(\U_s)_{s\in\N}$ of order~$n+1$.
\end{proposition}
\begin{proof}
We proceed as in the proof of Lemma \ref{lem_good_covers}, but when searching for a good cover~$\U_{s+1}$ that strongly refines~$\U_s$, we additionally test whether~$\U_{s+1}$ has order~$n+1$, which is also a~$\Sigma^0_1$ predicate. As~$\dim(X)\leq n$, it always exists so it can be effectively found.
\end{proof}

\begin{corollary}\label{cor_tree}
If~$X$ is zero-dimensional, then from any compact presentation of~$X$, one can compute a pruned tree~$T\subseteq 2^{<\omega}$ such that~$X$ is homeomorphic to~$[T]$.
\end{corollary}
\begin{proof}
As~$X$ is zero-dimensional, one can compute by Proposition \ref{prop_order} a strongly refining sequence of good open covers~$(\U_s)_{s\in\N}$ of order~$1$. Each~$\U_s$ is therefore made of clopen sets. All these clopen sets ordered by inclusion form a finitely-branching pruned tree, in which the number of nodes at each level can be computed. By standard arguments, it can be effectively converted into a binary tree.
\end{proof}

\section{Cantor-Bendixson derivative}\label{sec:Cantor-Bendixson}

Let~$X$ be a topological space. The {\em Cantor-Bendixson derivative} of $X$ is the subspace~$X'$ of all non-isolated points of $X$. We discuss the computability of~$X'$ in comparison with the computability of~$X$. This problem has been investigated in the context of reverse mathematics in \cite{GrMo08}.

In the next result, it is important to note that in order to produce a compact presentation of~$X'$, we only use a Polish presentation of~$X$.

\begin{lemma}\label{lem:CB-derivative}
Let ${X}$ be a compact Polish space.

If~$X\subseteq\hil$ is~$\d$-computably overt then~$X'$ is~$\Pi^0_1(\d')$ and~$\d''$-computably overt.

Therefore, if~$X$ has a~$\d$-computable Polish presentation, then its Cantor-Bendixon derivative~$X'$ has a~$\d''$-computable compact presentation.
\end{lemma}

\begin{proof}
Assume that~$X\subseteq Q$ is~$\d$-computably overt, and let~$(x_i)_{i\in\N}$ be a~$\d$-computable dense sequence in~$X$.

Let $B=B(s,r)$ be a rational open ball in $\hil$, and for~$k\in\N$ let~$B^k=B(s,r-2^{-k})$. We claim that $B$ intersects~$X'$ if and only if there is $k$ such that $B^k\cap X$ contains infinitely many points. For the forward direction, choose $x\in B\cap X'$. One has~$x\in B^k$ for sufficiently large~$k$. Since $B^k$ is open and~$x$ is not isolated in~$X$,~$B^k$ contains infinitely many points. For the backward direction, if~$B^k\cap X$ contains infinitely many points then its closure, which is contained in~$B\cap X$, contains a non-isolated point, therefore~$B$ intersects~$X'$. By this claim, the property $B\cap{X}'\not=\emptyset$ is equivalent to
\begin{equation*}
\exists k,\forall n,\exists x_{i_1},\ldots,x_{i_n}\text{ that are pairwise distinct and belong to~$B^k$},
\end{equation*}
which is~$\Sigma^0_3$ relative to~$\d$, or equivalently $\Sigma^0_1$ relative to~$\d''$. This means that $X'$ is $\d''$-computably overt.

Next, let~$A=\{(i,k):\forall j,x_j=x_i\text{ or }d_Q(x_i,x_j)\geq 2^{-k}\}$. The set~$A$ is~$\Pi^0_1$ relative to~$\d$, hence~$\d'$-computable. One can easily see that~$x$ is isolated in~${X}$ if and only if~$x\in B(x_i,2^{-k})$ for some~$(i,k)\in A$. Thus, the set of isolated points is a~$\Sigma^0_1(\d')$ subset of~${X}$; hence~${X}'$ is a~$\Pi^0_1(\d')$ subset of~${X}$, which is itself~$\Pi^0_1(\d')$, so~$X'$ is~$\Pi^0_1(\d')$.

Therefore,~$X'$ is~$\d''$-computably compact.
\end{proof}
%
%

\subsection{Stone duality}\label{boolean}

Here we show that spectra of compact zero-dimensional spaces are closely related to spectra of Boolean algebras. This follows from an effectivization of Stone duality in \cite{OdSe89}.

Let $\mathbf{B}$ be the category formed by the Boolean algebras as
objects and the $\{\vee,\wedge,\bar{}\;,0,1\}$-homomorphisms as
morphisms. Recall that a {\em Stone space}  is a compact topological
space $X$ such that for any distinct $x,y\in X$ there is a
clopen set $U$ with $x\in U\not\ni y$ (i.e., zero-dimensional and $T_1$).
Let $\mathbf{S}$ be the
category formed by the Stone spaces as objects and the continuous
mappings as morphisms.

The {\em Stone duality} states the dual
equivalence between the categories $\mathbf{B}$ and $\mathbf{S}$.
More explicitly, the Stone space $s(B)$  corresponding to a
given Boolean algebra $B$ is formed by the set of prime
filters of $B$ with the base of open (in fact, clopen)
sets consisting of the sets $\{F\in s(B)\mid a\in F\}$,
$a\in B$. (Note that one could equivalently take ideals in place
of filters.) Conversely, the Boolean algebra
$b(X)$ corresponding to a given Stone space $X$ is formed by the
set of clopen sets (with the usual set-theoretic operations). By
Stone duality, any Boolean algebra $B$ is canonically
isomorphic to the Boolean algebra $b(s(B))$ (the
isomorphism $f:B\rightarrow b(s(B))$ is defined
by $f(a)=\{F\in s(B)\mid a\in F\}$), and  any Stone
space $X$ is canonically homeomorphic to the space $s(b(X))$. 

Restricting the Stone duality to the countable Boolean algebras, we obtain their duality with the class~$\CP_0$ of compact zero-dimensional countably based spaces, and in fact with the compact subspaces of the Cantor space $2^\omega$. As the nonempty closed subsets of  $2^\om$ coincide with the sets $[T]$ of infinite paths through a pruned tree $T\subseteq2^\omega$, we obtain a close relation between such subspaces and countable Boolean algebras.

\begin{fact}\label{fact:Boolean-algebra}
\begin{enumerate}
\item A Boolean algebra has a $\d$-c.e.~(resp. $\d$-co-c.e., $\d$-computable) copy if and only if it is isomorphic to the Boolean algebra of clopen subsets of $[T]$ for some $\d$-co-c.e.~(resp. $\d$-c.e., $\d$-computable) pruned tree $T$.
\item Every $\d$-co-c.e.~Boolean algebra is isomorphic to a $\d$-computable Boolean algebra.
\item There is a $\d$-c.e.~Boolean algebra which is not isomorphic to a $\d$-computable Boolean algebra. 
\end{enumerate}
\end{fact}

\begin{proof}
The first item follows from \cite[Lemma 3]{OdSe89}; see also \cite{Se88}.
The second item follows from \cite{OdSe89}.
The third item follows from \cite{Fe70}.
\end{proof}

As already noticed by Harrison-Trainer, Melnikov, and Ng \cite{HTMN20}, one can use Stone duality to show several results on degree spectra of zero-dimensional compacta.
For instance, Stone duality can be used to show the following:

\begin{fact}[see \cite{HTMN20}]
\begin{enumerate}
\item There exists~$X\in\CP_0$ which has a $\zero'$-computable Polish presentation, but no computable Polish presentation, 
\item If~$X\in\CP_0$ has a low$_4$ Polish presentation, then it has a computable Polish presentation.
\end{enumerate}
\end{fact}

For~$X\in\CP_0$, Corollary \ref{cor_tree} implies that~$X$ has a~$\d$-computable compact presentation if and only if there exists a~$\d$-computable pruned tree~$T\subseteq 2^{<\omega}$ such that~$X$ is homeomorphic to~$[T]$. Therefore,
\begin{corollary}\label{cor_ba_zero}
The degree spectra of countable boolean algebras coincide with the compact degree spectra of zero-dimensional compact Polish spaces.
\end{corollary}

The Stone dual of the Cantor-Bendixon derivative is known as the Fr\'echet derivative $B'$ of a Boolean algebra $B$ which is the quotient of $B$ by the ideal generated by atoms (minimal non-zero elements). Since the isolated points $x$ of the space $s(B)$ (realized as $[F]$ above) are precisely the atoms $[\tau]\cap[F]$ for suitable prefix $\tau\sqsubseteq x$, we obtain the following.

\begin{proposition}\label{deriv}
For any countable Boolean algebra $B$, $s(B')$ is homeomorphic to $s(B)'$.
\end{proposition}

Precise complexity estimations for the Frechet derivative were obtained in \cite{OdSe89}: for any Turing degree~$\d$, a countable Boolean algebra $C$ is $\d''$-computably presentable iff $C$ is isomorphic to $B'$ for some  $\d$-computable Boolean algebra $B$, and there is a $\d$-computable Boolean algebra $B$ such that $B'$ is not $\d'$-computably presentable. 
The iterated version is also known for any $n>0$: a countable Boolean algebra $C$ is $\d^{(2n)}$-computably presentable iff $C$ is isomorphic to the $n$th derivative $B^{(n)}$ for some $\d$-computable Boolean algebra $B$, and there is a $\d$-computable Boolean algebra $B$ such that the $n$th derivative $B^{(n)}$ is not $\d^{(2n-1)}$-computably presentable.

These results have an immediate consequence in terms of compact degree spectra of spaces in~$\CP_0$, by Corollary \ref{cor_ba_zero} and Proposition \ref{deriv}.

\begin{theorem}\label{thm:Cantor-Bendixson-CP0}
Let~$\d$ be any Turing degree. For any $n\geq 1$, a space $Y\in\mathbf{CP}_0$ has a $\d^{(2n)}$-computable compact presentation if and only if $Y$ is homeomorphic to the $n$th derivative $X^{(n)}$ for some $\d$-computable compact $X\in\mathbf{CP}_0$, and there is a $\d$-computable compact $X\in\mathbf{CP}_0$ such that the $n$th derivative $X^{(n)}$ does not have a $\d^{(2n-1)}$-computable compact presentation.
\end{theorem}

 It implies in particular that Lemma \ref{lem:CB-derivative} is almost optimal.
 
\subsection{Countable spaces}

We next note that countable topological spaces cannot be used to construct nontrivial degree spectra inside the hyperarithmetical hierarchy.

Mazurkiewicz-Sierpi\'nski's theorem states that every countable compact Polish space is homeomorphic to the ordinal space $\om^\alpha\cdot n+1$ for some $\alpha<\om_1$ and $n\in\om$, endowed with the order topology.

Let $\om_1^x$ be the least ordinal which is not computable in $x$.

\begin{proposition}
For any countable ordinal $\alpha$, the compact and Polish degree spectrum of the space~$\omega^\alpha+1$ are both $\{x:\alpha<\om_1^x\}$.
\end{proposition}

\begin{proof}
We first show that if $\alpha$ is $\x$-computable, then $\omega^\alpha+1$ has an $\x$-computable compact presentation.

Given an~$\x$-computable well-ordering~$\preceq$ of~$\N$ of order type~$\omega^\alpha+1$, we build a copy~$X$ of~$\omega^\alpha+1$ in~$[0,1]$ as follows. We can assume w.l.o.g.~that~$0$ is minimal and~$1$ is maximal for~$\preceq$. We embed~$(\N,\preceq)$ in~$([0,1],\leq)$ by inductively defining a rational number~$x_n\in [0,1]$ for each~$n\in\N$. We need some care to make sure that we will obtain a topological embedding, and that the copy~$X$ will be~$\x$-computably compact. We start with~$x_0=0$ and~$x_1=1$. Let~$n\geq 2$ and assume that~$x_0,\ldots,x_{n-1}$ have been defined. Let~$i$ and~$j$ be respectively  the predecessor and the successor of~$n$ in~$(\{0,\ldots,n\},\preceq)$, and define~$x_{n}=2^{-n}x_i+(1-2^{-n})x_j$. Note that~$x_j-x_n\leq 2^{-n}$. Let~$X=\{x_n:n\in\N\}$ with the topology inherited from~$\R$.

The function~$n\mapsto x_n$ is by construction an order embedding of~$(\N,\preceq)$ in~$([0,1],\leq)$. We need to show that it is a topological embedding, i.e.~that it is continuous. It is sufficient to check that if~$n=\sup_\preceq \{m:m\prec n\}$, then~$x_n=\sup_\leq\{x_m:m\prec n\}$. There exist infinitely many~$m$'s such that~$m\geq n$,~$m\prec n$ and~$m$ is the predecessor of~$n$ in~$(\{0,\ldots,m\},\preceq)$. When defining~$x_m$ for any such~$m$, the successor of~$m$ is~$j=n$ so~$x_n-x_m\leq 2^{-m}$, which is arbitrarily small so~$x_n=\sup_\leq\{m:m\prec n\}$ as wanted.

 
Therefore, the set~$X=\{x_n:n\in\N\}$ is a copy of~$\omega^\alpha+1$. The sequence~$(x_n)_{n\in\N}$ is an~$\x$-computable sequence of rational numbers. Let~$X_n=\{x_0,\ldots,x_n\}\subseteq X$ and observe by construction that~$d_H(X_n,X_{n-1})\leq 2^{-n}$, so~$d_H(X_n,X)\leq 2^{-n}$. As the compact sets~$X_n$ are uniformly~$\x$-computable and converge fast to~$X$,~$X$ is~$\x$-computable as well.

Conversely, if $\om^\alpha+1$ has an $\x$-computable Polish presentation, then by Lemma \ref{lem:Polish-to-compact}, it has an $\x'$-computable compact presentation. In particular, there is a countable $\Pi^0_1(\x')$ class $P\subseteq 2^\omega$ which is homeomorphic to $\om^\alpha+1$. Since the Cantor-Bendixson rank of $\om^\alpha+1$ is $\alpha$, and the Cantor-Bendixson rank is a topological invariant, the rank of $P$ is also $\alpha$. However, as noted by Kreisel, the Spector boundedness principle implies that the rank of a countable $\Pi^0_1(\x')$ class must be $\x'$-computable; see also \cite[Section 4]{Ce99}. As an $\x$-hyperarithmetic ordinal is always $\x$-computable, this implies that $\alpha<\om_1^{\bf x}$.
\end{proof}

The same arguments show that the compact and Polish degree spectra of~$\omega^\alpha\cdot n+1$ are both $\{x:\alpha<\om_1^x\}$. Therefore, if a countable Polish space has an hyperarithmetical compact (or even Polish) presentation, then it has a computable compact presentation.

For more details, see also \cite{AK00,GrMo08}.

\section{Computable aspects of \v{C}ech homology groups}\label{sec_homo}

\subsection{Background on \v{C}ech homology groups}
All the details about the definitions and results mentioned in this section can be found in the textbook of Hurewicz and Wallman \cite{HurWal41}. We will refer to the corresponding sections in that reference.

\paragraph{Simplicial complexes.} We refer to \cite[Section V.9]{HurWal41}.
Let~$V$ be a finite set, whose elements are called the \emph{vertices}. A \emph{simplicial complex}~$K$ on~$V$ is a collection of subsets of~$V$ such that if~$A\subseteq B\subseteq V$ and~$B\in K$, then~$A\in K$. An element~$A\in K$ is called a \emph{simplex}, and an~$n$\emph{-simplex} if the cardinality of~$A$ is~$n+1$.

To a simplicial complex~$K$ corresponds a topological space~$X_K$, called its \emph{geometric realization}. Assuming that~$V=\{0,\ldots,k\}$, let~$\Delta_V=\{(x_0,\ldots,x_k):\forall i, x_i\in [0,1]\text{ and }\sum_i x_i=1\}$. For~$x=(x_0,\ldots,x_k)\in\Delta_V$, let~$\supp(x)=\{i:x_i>0\}$. The geometric realization of~$K$ is the subset of~$\Delta_V$ defined by
\begin{equation*}
X_K:=\{x\in \Delta_V:\supp(x)\in K\}.
\end{equation*}

Let~$K,L$ be simplicial complexes on vertices~$V_K$ and~$V_L$ respectively. A \emph{simplicial map}, written~$f:K\to L$, is a map~$f:V_K\to V_L$ such that if~$A\in K$ then~$f(A)=\{f(v):v\in A\}\in L$. In other words, a simplicial map sends each simplex of~$K$ to a simplex of~$L$. Note that if~$A$ is an~$n$-simplex of~$K$, then~$f(A)$ is a~$p$-simplex of~$L$ for some~$p\leq n$, and~$f(A)$ is an~$n$-simplex if and only if the images of the vertices of~$A$ under~$f$ are pairwise distinct.


\paragraph{Simplicial homology groups.} We refer to \cite[Sections VIII.1, VIII.3]{HurWal41}. To a simplicial complex~$K$ are associated its homology groups~$H_n(K;G)$ where~$G$ is an abelian group, called the coefficient group. Here we will only consider the coefficient group~$G=\zd:=\ZZ/2\ZZ$.

Let~$K$ be a simplicial complex on a finite set~$V$ of vertices. An~$n$\emph{-chain} of~$K$ is a set of~$n$-simplices of~$K$, or equivalently a formal sum of~$n$-simplices with coefficients in~$\zd$. Let~$C_n(K)$ be the finite abelian group of~$n$-chains, endowed with the sum operation. The \emph{boundary operator} is the group homomorphism~$\partial_n:C_{n}(K)\to C_{n-1}(K)$ that sends an~$n$-simplex~$\Delta=\{v_0,\ldots,v_{n}\}$ to~$\sum_{i=0}^{n}\Delta\setminus \{v_i\}$.

An~$n$\emph{-cycle} is an~$n$-chain whose boundary is~$0$, i.e.~is an element of~$\ker\partial_{n}$. An~$n$\emph{-boundary} is an~$n$-chain which is the boundary of some~$(n+1)$-chain, i.e.~is an element of~$\im \partial_{n+1}$. One has~$\partial_{n}\circ \partial_{n+1}=0$, in other words~$\im\partial_{n+1}\subseteq\ker\partial_{n}$, i.e.~every~$n$-boundary is an~$n$-cycle. The $n$th \emph{homology group} of~$K$ is the finite abelian group~$H_n(K;\zd)=\ker\partial_{n}/\im \partial_{n+1}$.

Let~$K,L$ be simplicial complexes and~$f:K\to L$ a simplicial map. For each~$n\in\N$,~$f$ induces group homomorphisms~$f_\#:C_n(K)\to C_n(L)$ defined as follows. If~$A\in K$ is an~$n$-simplex, then
\begin{equation*}
f_\#(A)=\begin{cases}
f(A)&\text{if~$f(A)$ is an~$n$-simplex,}\\
0&\text{otherwise.}
\end{cases}
\end{equation*}
Said differently,~$f_\#(A)=f(A)$ if the images of the vertices of~$A$ are pairwise distinct. The homomorphism~$f_\#$ commutes with~$\partial_{n+1}$, i.e.~$f_\#\circ \partial_{n+1}=\partial_{n+1}\circ f_\#$, which implies that~$f_\#$ induces a group homomorphism~$f_*:H_n(K;\zd)\to H_n(L;\zd)$.

\paragraph{\v{C}ech homology groups.} We refer to \cite[Section VIII.4]{HurWal41}. There are several ways to define homology groups of more general topological spaces. We will use the \v{C}ech definition, which approximates the space by simplicial complexes.

Let~$X$ be a compact Polish space. To a finite open cover~$\U=(U_0,\ldots,U_k)$ of~$X$ is associated its \emph{nerve}~$K_\U$, which is a simplicial complex. Its set of vertices is~$V_\U=\{0,\ldots,k\}$, and a simplex~$A\subseteq V_\U$ belongs to~$K_\U$ if and only if~$\bigcap_{i\in A}U_i\neq \emptyset$.

Let~$(\U_s)_{s\in\N}$ be a sequence of finite open covers of~$X$ such that~$\U_{s+1}$ refines~$\U_s$, i.e.~each open set of~$\U_{s+1}$ is contained in some open set of~$\U_s$, and such that the maximal diameter of elements of~$\U_s$ tends to~$0$ as~$s$ grows.

Let~$K_s$ be the nerve of~$\U_s$. As~$\U_{s+1}$ refines~$\U_s$, one can define a simplicial map~$f_{s}:K_{s+1}\to K_s$ by sending a vertex corresponding to~$U\in \U_{s+1}$ to a vertex corresponding to some~$V\in\U_s$ containing~$U$. There are several possible choices for~$f_s$, but in the end they will lead to the same result.

The $n$th \emph{\v{C}ech homology group}~$\check{H}_n(X;\zd)$ of~$X$ is defined as the inverse limit of the inverse system~$f_{s,*}:H_n(K_{s+1};\zd)\to H_n(K_{s};\zd)$. In other words,~$\check{H}_n(X;\zd)$ is the set of sequences~$(c_s)_{s\in\N}$ with~$c_s\in H_n(K_s;\zd)$ and~$c_s=f_{s,*}(c_{s+1})$, with component-wise addition. This group does not depend on the choice of open covers and simplicial maps.

If~$X$ is itself homeomorphic to the geometric realization of a simplicial complex, then the \v{C}ech homology group~$\check{H}_n(X;\zd)$ is isomorphic to the simplicial homology group~$H_n(X;\zd)$ \cite[Theorem VIII.4.E]{HurWal41}.

\subsection{Computable aspects}
Lupini, Melnikov, Nies \cite{LMN23} and Downey, Melnikov \cite{DM23} proved that the \v{C}ech cohomology groups~$\check{H}^n(X)$ of a compact Polish space~$X$ are computable, i.e.~can be presented as a computable set of relations over a countable set of generators, relative to any compact presentation of~$X$. We will use part of their argument to investigate the computability of the \v{C}ech homology groups~$\check{H}_n(X;\zd)$. However, these groups do not behave so well in terms of computability. First, they are not countable in general, so they have no countable set of generators and their computability does make immediate sense. Second, as we will show, their non-triviality is~$\Sigma^1_1$-complete when taking coefficients in~$\ZZ$, which contrasts with the non-triviality of the \v{C}ech cohomology groups with coefficients in~$\ZZ$, which is~$\Sigma^0_2$. However, we now show that when taking coefficients in a finite group such as~$\zd$, the non-triviality of the \v{C}ech homology groups is~$\Sigma^0_2$.

Given a combinatorial description of a finite simplicial complex~$K$, each group~$H_n(K;\zd)$ can be fully computed: it is a finite group whose cardinality and operation can be computed, for instance under the form of a Cayley table. Given a combinatorial description of~$K,L$ and of a simplicial map~$f:K\to L$, the group homomorphisms~$f^*_n:H_n(K;\zd)\to H_n(L;\zd)$ can be computed, as functions between finite sets. It makes the inverse limit~$\check{H}_n(X;\zd)$ a co-c.e.~profinite group, in the sense of La Roche \cite{LaRoche81} and Smith \cite{Smith81}, i.e.~the inverse limit of a computable inverse system of finite groups.

In the next statement, a finite group is represented in the strongest possible way by a Cayley table, which encodes the cardinality of the group and the group operation.   A group homomorphism is also represented in the strongest possible way by arranging its graph in a finite list.
\begin{proposition}\label{prop_nontrivial}
Given an inverse system of finite groups~$f_n:G_{n+1}\to G_{n}$, whether their inverse limit is non-trivial is~$\Sigma^0_2$.
\end{proposition}
\begin{proof}
Let~$G$ be the inverse limit. We claim that~$G$ is non-trivial if and only if there exist~$n\in\N$ and~$g\in G_n$ such that~$g\neq 0$ and~$g\in\im(f_n\circ\ldots\circ f_p)$ for all~$p\geq n$. This predicate is~$\Sigma^0_2$, because the groups are finite. In particular, the predicate~$g\in\im(f_n\circ \ldots\circ f_p)$ can be decided in finite time, by testing~$f_n\circ \ldots\circ f_p(g')=g$ for each~$g'$ in the finite set~$G_{p+1}$.

If~$G$ is non-trivial, then let~$(g_n)_{n\in\N}$ be a non-zero element of~$G$. As it is non-zero, there exists~$n$ such that~$g_n\neq 0$. Moreover, for every~$p\geq n$,~$g_n=f_n\circ\ldots\circ f_p(g_{p+1})$ belongs to~$\im(f_n\circ\ldots\circ f_p)$.

Conversely, assume the existence of~$n$ and~$g$ satisfying the conditions. As the groups are finite, by K\"onig's lemma, there exists a sequence~$(g_p)_{p\geq n}$ such that~$g_n=g$ and~$g_p=f_p(g_{p+1})$ for all~$p\geq n$. For~$m<n$, let~$g_m=f_m\circ \ldots \circ f_{n-1}(g)$. The sequence~$(g_n)_{n\in\N}$ is an element of~$G$, and is non-zero as its~$n$th coordinate is~$g\neq 0$.
\end{proof}


We now have all the ingredients to show that the non-triviality of \v{C}ech homology groups with coefficients in~$\zd$ has limited complexity.

\begin{corollary}\label{cor_nontrivial_homology}
Whether the~$n$th \v{C}ech homology group~$\check{H}_n(X;\zd)$ of a compact Polish space~$X$ is non-trivial is~$\Sigma^0_2$, uniformly in~$n$ and a compact presentation of~$X$.
\end{corollary}
\begin{proof}
Given a compact presentation of~$X$ and a good open cover~$\U=(U_0,\ldots,U_k)$ of~$X$, its nerve~$K_\U$ on the vertices~$V=\{0,\ldots,k\}$ can be computed, in the sense that one can decide whether a set~$J\subseteq V$ belongs to~$K_\U$. Indeed,~$J\in K_\U$ iff~$\bigcap_{j\in J}U_j\neq\emptyset$ iff~$\bigcap_{j\in J}\overline{U}_j\neq \emptyset$. The first inequality is~$\Sigma^0_1$ and the second inequality is~$\Pi^0_1$, so it is decidable.

If~$\U=(U_0,\ldots,U_k)$ and~$\V=(V_0,\ldots,V_l)$ are good open covers of~$X$ and~$\U$ strongly refines~$\V$, then one can compute a simplicial map~$f:K_\U\to K_\V$ sending each vertex~$u_i$ for~$K_\U$ to a vertex~$v_j$ of~$K_\V$ such that~$\overline{U}_i\subseteq V_j$. Indeed, the latter inclusion is~$\Sigma^0_1$, so for each~$i$ one can effectively find a suitable~$j$.

Now, from a compact presentation of~$X$, one can compute a strong refining sequence of good open covers~$(\U_s)_{s\in\N}$ by Lemma \ref{lem_good_covers}. One can then compute the nerves~$K_s$ of~$\U_s$ as well as simplicial maps~$f_s:K_{s+1}\to K_s$. One can compute the simplicial homology groups~$H_n(K_s)$ and the homomorphisms~$f_{s,*}:H_n(K_{s+1};\zd)\to H_n(K_s;\zd)$. Their inverse limit is~$\check{H}_n(X;\zd)$, and whether it is non-trivial is~$\Sigma^0_2$ by Proposition \ref{prop_nontrivial}.
\end{proof}

\paragraph{Comparison between \v{C}ech homology and \v{C}ech cohomology.}
Lupini, Melnikov and Nies \cite{LMN23} and Downey and Melnikov \cite{DM23} proved that the \v{C}ech \emph{cohomology} groups~$\check{H}^n(X)$ of a compact Polish space~$X$ can be computed from a compact presentation of~$X$ (we do not write the coefficient group, which is~$\ZZ$). We will neither define nor use cohomology groups, but let us mention the important fact that for each~$n$, the simplicial cohomology groups~$H^n(K_s)$ and the induced homomorphisms~$f_s^*:H^n(K_{s})\to H^{n}(K_{s+1})$ form a \emph{direct} system, which is essential to show that their direct limits, the \v{C}ech cohomology groups~$\check{H}^n(X)$, are computable. As we saw, the simplicial homology groups~$H_n(K_s;\zd)$ form an \emph{inverse} system, which makes the computation of their inverse limits, the \v{C}ech \emph{homology} groups~$\check{H}_n(X;\zd)$, much more difficult.



\subsection{Influence of the coefficient group}
 Let us finally discuss the role of the coefficient group in Corollary \ref{cor_nontrivial_homology}. The same proof holds for any finite coefficient group~$G$. Moreover, the computability of the \v{C}ech cohomology groups proven in \cite{LMN23,DM23} implies that Corollary \ref{cor_nontrivial_homology} holds for the circle group~$G=\R/\ZZ$. Indeed, there is a duality between \v{C}ech cohomology groups and \v{C}ech homology groups:~$\check{H}_n(X;\R/\ZZ)$ is the character group of~$\check{H}^n(X)$, i.e.~is the group of homomorphisms from~$\check{H}^n(X)$ to the circle group~$\R/\ZZ$ (\cite[Theorem VIII.4.G, p.~137]{HurWal41}). Therefore, the non-triviality of~$\check{H}_n(X;\R/\ZZ)$ is equivalent to the non-triviality of~$\check{H}^n(X)$, which is~$\Sigma^0_2$ in~$X$, as~$\check{H}^n(X)$ is computable in~$X$.

We now show that Corollary \ref{cor_nontrivial_homology} does not hold for~$G=\ZZ$, in which case the non-triviality of~$\check{H}_n(X;\ZZ)$ is~$\Sigma^1_1$-complete. The reason is that the groups of the inverse system are no more finite, so the existence of a path in the inverse system is similar to the existence of an infinite path in a tree.
\begin{theorem}
Let~$n\geq 1$. Given a compact presentation of a compact Polish space~$X$, the non-triviality of~$\check{H}_n(X;\ZZ)$ is~$\Sigma^1_1$-complete.
\end{theorem}

\begin{proof}
The arguments in \cite{Melnikov21,LMN23,DM23} imply the computability of the inverse system of simplicial groups~$(H_n(K_s;\ZZ),f_{s,*})$: these groups are \emph{strongly completely decomposable} (strictly speaking, their results are about cohomology groups~$H^n(K_s)$ but apply equally to the homology groups~$H_n(K_s;\ZZ)$). Therefore, their inverse limit is non-trivial if and only if there exists a sequence~$(x_s)_{s\in\N}$ with~$x_s\in H_n(K_s;\ZZ)$,~$x_{s}=f_{s,*}(x_{s+1})$ for all~$s$ and~$x_s\neq 0$ for some~$s$. This condition is~$\Sigma^1_1$. We now show the completeness of the problem.

The set of trees~$T\subseteq 2^{<\omega}$ containing an infinite path with infinitely many~$1$'s is~$\Sigma^1_1$-complete \cite[Proposition 25.2]{Kechris95}. We will reduce this problem to the non-triviliaty of~$\check{H}_1(X;\ZZ)$, and explain the changes to be made for higher-dimensional homology groups. We will build an inverse system of groups and then build an inverse limit of compact spaces whose \v{C}ech homology group is the inverse limit of the groups. Note the similarity of this strategy with \cite{Melnikov21} which builds an inverse limit of compact spaces whose first \v{C}ech cohomology group is a given \emph{direct} limit of groups.
 

We first define an inverse system of groups~$(G_n,f_n)$, i.e.~a sequence~$(G_n)_{n\in\N}$ of groups together with homomorphisms~$f_n:G_{n+1}\to G_n$. We will then associate to each tree~$T$ a subsystem~$(L_n,f_n)$, i.e.~a sequence of subgroups~$L_n\subseteq G_n$ satisfying~$f_n(L_{n+1})\subseteq L_n$.

Let~$\epsilon\in 2^{<\omega}$ be the empty string,~$S=\{\epsilon\}\cup \{1w:w\in 2^{<\omega}\}$ be the set of strings that do not start with~$0$ and~$S_n=\{w\in S:|w|=n\}$. For each~$w\in S$, let~$k(w)=2^{|w|_1}$ where~$|w|_1$ is the number of occurrences of~$1$ in~$w$. Note that~$k(\epsilon)=1$.

For~$n\in\N$, we define
\begin{equation*}
G_n=\bigoplus_{w\in S_n}\ZZ^{k(w)}.
\end{equation*}
Note that~$G_0=\ZZ$ and~$G_1=\ZZ^2$. We call each~$\ZZ^{k(w)}$ a component of~$G_n$. If~$n\geq 1$ and~$w\in S_n$, then~$|w|_1\geq 1$ so~$k(w)$ is even: we will implicitly split each component~$\ZZ^{k(w)}$ into~$\ZZ^{k(w)/2}\oplus \ZZ^{k(w)/2}$, expressing an element of~$\ZZ^{k(w)}$ as a pair.

We define~$f_n:G_{n+1}\to G_n$ on each component of~$G_{n+1}$. For each~$w\in S_n$, one has~$k(w0)=k(w)$ and~$k(w1)=2k(w)$ and we let
\begin{align*}
f_n:\ZZ^{k(w0)}&\to \ZZ^{k(w)}\\
(x,y)&\mapsto (2x,3y)
\end{align*}
and
\begin{align*}
f_n:\ZZ^{k(w1)}&\to\ZZ^{k(w)}\\
(x,y)&\mapsto x+y.
\end{align*}

An element~$x\in\varprojlim (G_n,f_n)$ is described as a family~$(x(w))_{w\in S}$ where~$x(w)\in\ZZ^{k(w)}$ for each~$w\in S$, such that~$x(w)=f_{|w|}(x(w0))+f_{|w|}(x(w1))$. We define the support of~$x$ as~$\supp(x)=\{w\in S:x(w)\neq 0\}$. Note that~$x=0$ if and only if~$\supp(x)=S$.

We relate the supports of non-zero elements with the infinite binary sequences containing infinitely many~$1$'s. For an infinite sequence~$q\in 2^\omega$, let~$q\res{n}$ be the prefix of~$q$ of length~$n$.
\begin{claim}\label{claim_supp1}
Let~$x\in\varprojlim (G_n,f_n)$. One has~$x\neq 0$ iff there exists~$q\in 2^\omega$ starting with~$1$ and containing infinitely many~$1$'s, and such that~$x(q\res{n})\neq 0$ for all sufficiently large~$n$.
\end{claim}
\begin{proof}
Of course if such a~$q$ exists, then~$x\neq 0$. Assume that~$x\neq 0$ and let~$w\in S$ be such that~$x(w)\neq 0$. We build~$q$ extending~$w$. Let~$n\geq |w|$ and assume by induction that~$v:=q\res{n}$ has been already defined. As~$x(v)=f_{n}(x(v0))+f_{n}(x(v1))$, one has~$x(v0)\neq 0$ or~$x(v1)\neq 0$. We choose the next bit of~$q$ so that~$x(q\res{n+1})\neq 0$: let~$q(n)=1$ if~$x(v1)\neq 0$, and~$q(n)=0$ otherwise. One has~$x(q\res{n})\neq 0$ for all~$n\geq |w|$ by construction. Assume for a contradiction that~$q$ contains only finitely many~$1$'s. Let~$v$ be a prefix of~$q$ extending~$w$ and such that~$q=v0^\omega$. Note that for~$n\geq |v|$,~$k(q\res{n})$ is constantly~$k(v)$ and~$x(q\res{n})\in \ZZ^{k(v)}$. For such~$n$, one must have~$x(q\res{n}1)=0$ because otherwise we would have chosen~$q(n)=1$ in the construction of~$q$. As a result, one has~$x(q\res{n})=f_n(x(q\res{n}0))=f_n(x(q\res{n+1}))$ for all~$n\geq |v|$. It means that each coordinate of~$x(q\res{v})\in \ZZ^{k(v)}$ is a multiple of all~$2^i$ for all~$i\geq 0$, or a multiple of~$3^i$ for all~$i\geq 0$, which implies that~$x(q\res{v})=0$. It is a contradiction because~$v$ extends~$w$. Therefore, we have proved that~$q$ contains infinitely many~$1$'s.
\end{proof}

Conversely,
\begin{claim}\label{claim_supp2}
If~$q\in2^\N$ starts with~$1$ and contains infinitely many~$1$'s, then there exists~$x\in\varprojlim (G_n,f_n)$ such that~$\supp(x)$ is the set of prefixes of~$q$.
\end{claim}
\begin{proof}
Let~$(n_i)_{i\in\N}$ be the increasing enumeration of the positions~$n$ such that~$q(n)=1$, and observe that~$n_0=0$. We build~$x=(x(w))_{w\in S}$. Let~$x(\epsilon)=1\in\ZZ$. Let~$i\in\N$ and assume by induction that~$x(q\res{n_i})$ has been defined. The function~$f_{n_i}\circ\ldots\circ f_{n_{i+1}-1}:\ZZ^{k(q\res{n_{i+1}})}\to\ZZ^{k(q\res{n_i})}$ maps~$(a,b)$ to~$2^ja+3^jb$, where~$j=n_{i+1}-n_i-1$. As~$2^j$ and~$3^j$ are coprime, there exists~$(a,b)$ such that~$2^ja+3^jb=x(q\res{n_i})$, let~$x(q\res{n_{i+1}})=(a,b)$. For~$n_i<n<n_{i+1}$, naturally define~$x(q\res{n})=f_n\circ\ldots f_{n_{i+1}-1}(a,b)$. Finally, let~$x(w)=0$ if~$w$ is not a prefix of~$q$. By construction,~$x$ belongs to~$\varprojlim (G_n,f_n)$ and~$\supp(x)$ is the set of prefixes of~$q$.
\end{proof}

Let~$T\subseteq 2^{<\omega}$ be a tree. We can assume w.l.o.g.~that~$T\subseteq S$, replacing~$T$ with~$\{\epsilon\}\cup\{1w:x\in T\}$ if necessary. We define an inverse subsystem~$(L_n,f_n)$ of~$(G_n,f_n)$. Let~$T_n=\{w\in T:|w|=n\}$ and~$L_n=\bigoplus_{w\in T_n}\ZZ^{k(w)}$. As~$T$ is a tree, one has~$f_n(L_{n+1})\subseteq L_n$ so~$(L_n,f_n)$ is an inverse system which is a subsystem of~$(G_n,f_n)$. Its inverse limit consists of the elements of the inverse limit of~$(G_n,f_n)$ whose support is contained in~$T$.

Therefore,~$\varprojlim (L_n,f_n)$ contains a non-zero element iff~$\varprojlim (G_n,f_n)$ contains a non-zero element whose support is contained in~$T$ iff~$T$ contains an infinite path with infinitely many~$1$'s by Claims \ref{claim_supp1} and \ref{claim_supp2}.

We now define a compact space~$X_T$ such that~$\check{H}_1(X_T;\ZZ)$ is isomorphic to~$\varprojlim (L_n,f_n)$. Let~$\T_n=\Sph_1^n$ be the $n$-dimensional torus, which is the product of~$n$ circles.

Observe that~$L_n=\ZZ^{l_n}$ for some~$l_n$. We then define a space~$X_T$ as the inverse limit of~$X_n=\T_{l_n}$, with continuous surjective functions~$F_n:X_{n+1}\to X_n$ defined as follows. Each homomorphism~$f_n:\ZZ^{l_{n+1}}\to \ZZ^{l_n}$ can be decomposed as a direct sum of the homomorphisms
\begin{align*}
g_2:\ZZ&\to \ZZ&g_3:\ZZ&\to\ZZ&g_+:\ZZ^2&\to\ZZ\\
x&\mapsto 2x&x&\mapsto 3x&(x,y)&\mapsto x+y.
\end{align*}
We define~$F_n$ as the corresponding product of
\begin{align*}
G_2:\Sph_1&\to\Sph_1&G_3:\Sph_1&\to\Sph_1&G_+:\Sph_1\times\Sph_1&\to\Sph_1\\
x&\mapsto 2x&x&\mapsto 3x&(x,y)&\mapsto x+y,
\end{align*}
where the circle~$\Sph_1$ is seen as an additive group. Let then~$X_T=\varprojlim (X_n,F_n)$.

The homology group of~$X_n=\T_{l_n}$ is~$H_1(X_n;\ZZ)=\ZZ^{l_n}=L_n$ and the homomorphism~$F_{n,*}:H_1(X_{n+1};\ZZ)\to H_1(X_n;\ZZ)$ induced by~$F_n$ is~$f_n$. It is a classical result that the construction of \v{C}ech homology groups of compact spaces with coefficients in~$\ZZ$ is \emph{continuous}, i.e.~commutes with inverse limits:
\begin{equation*}
\check{H}_1(\varprojlim(X_n,F_n);\ZZ)=\varprojlim (H_1(X_n;\ZZ),f_n)=\varprojlim (L_n,f_n).
\end{equation*}
This result can be found in \cite[Theorem X.3.1, p.261]{EilenbergS52}, where the statement holds for coefficient groups in the class~$\mathcal{G}_R$ of modules over a ring~$R$, and~$\ZZ$ is indeed a module over the ring~$\ZZ$. Therefore,~$\check{H}_1(X_T;\ZZ)=\varprojlim (L_n,f_n)$ is non-trivial if and only if~$T$ contains a path with infinitely many~$1$'s.

We finally need to check that the construction of~$X_T$ is effective in~$T$.

\begin{claim}
A compact presentation of the space~$X_T$ can be uniformly computed from~$T$.
\end{claim}
\begin{proof}
It is essentially the same argument as \cite[Remark 2.8]{Melnikov21}. For simplicity, we assume that~$T$ is computable, the general case holds relative to and uniformly in~$T$. First, the spaces~$X_n=\T_{l_n}$ have uniformly computable compact presentations, and so does their product~$P=\prod_n X_n$. The functions~$F_n:X_{n+1}\to X_n$ are uniformly computable so~$X_T$ is a~$\Pi^0_1$-subspace of~$P$. We need to show that~$X_T$ is computably overt in~$P$. It is important to note that the functions~$F_n:X_{n+1}\to X_n$ are surjective. For~$n\leq p$, let~$F_{n,p}=F_n\circ F_{n+1}\circ \ldots \circ F_{p-1}:X_p\to X_n$, with~$F_{p,p}=\mathrm{id}$. Let~$B=B_0\times \ldots\times B_p\times \prod_{q>p}X_q$ be a basic ball of~$P$. The open set~$B$ intersects~$X$ iff the set~$\bigcap_{n\leq p}F_{n,p}^{-1}(B_n)$ is non-empty, which is a c.e~condition. Indeed, if that set contains an element~$x_p$, then let~$x_n=F_{n,p}(x_p)$ for~$n<p$ and inductively choose for each~$q>p$ some~$x_q\in X_q$ such that~$F_{q-1}(x_q)=x_{q-1}$, which is possible as each~$F_{q-1}$ is surjective. The sequence~$(x_n)_{n\in\N}$ belongs to~$X\cap B$.
\end{proof}

Therefore, we have proved that the non-triviality of~$\check{H}_1(X;\ZZ)$ is~$\Sigma^1_1$-complete. We can achieve the same result for~$\check{H}_k(X;\ZZ)$ for any~$k\geq 1$. Given a binary tree~$T$, we build the same inverse system of groups~$(L_n,f_n)$, with~$L_n=\ZZ^{l_n}$. Let~$\Sph_k$ be the~$k$-dimensional sphere and note that~$H_k(\Sph_k;\ZZ)\cong\ZZ$. In the previous construction, we replace the circle by~$\Sph_k$. However, the equality~$H_k(\Sph_k^{l_n};\ZZ)\cong \ZZ^{l_n}$ only holds for~$k=1$, so we need to replace the product by the wedge sum (the previous argument would also work with the wedge sum). Let~$s\in\Sph_k$ be a distinguished point, making~$(\Sph_k,s)$ a pointed space, and let~$X_n=\bigvee_{i<l_n} (\Sph_k,s)$. One has~$H_k(X_n;\ZZ)=\ZZ^{l_n}$ \cite[Corollary 2.25]{Hatcher02}. 

We define continuous surjective functions~$F_n:X_{n+1}\to X_n$ in a similar way, by combining the following functions. Let~$G_2,G_3:(\Sph_k,s)\to(\Sph_k,s)$ be continuous functions of degrees~$2$ and~$3$ respectively, i.e.~inducing the maps~$g_2,g_3:H_k(\Sph_k;\ZZ)\to H_k(\Sph_k;\ZZ)$ respectively. Let~$G_+:(\Sph_k,s)\vee(\Sph_k,s)\to(\Sph_k,s)$ be defined as the identity on each copy of~$\Sph_k$. One has~$H_k(\Sph_k\vee\Sph_k)\cong H_k(\Sph_k)\oplus H_k(\Sph_k)\cong \ZZ^2$ and~$G_+$ induces~$g_+:\ZZ^2\to \ZZ$. We finally define~$X=\varprojlim (X_n,F_n)$. Again by continuity of \v{C}ech homology groups, one has~$\check{H}_k(X;\ZZ)=\varprojlim (H_k(X_n;\ZZ),f_n)=\varprojlim (L_n,f_n)$, which is non-trivial iff~$T$ contains a path with infinitely many~$1$'s.
\end{proof}

\section{Degree spectra of compact Polish spaces}\label{sec_spectra}
In this section, we realize certain sets of degrees as spectra of compact Polish spaces. We encode a set of natural numbers into the dimensions of spheres of a space. For each~$n\in\N$, the~$n$-dimensional sphere can be realized as~$\Sph^n=\{(x_0,\ldots,x_n)\in\R^{n+1}:x_0^2+\ldots+x_d^2=1\}$. To a set~$A\subseteq\N$ we associate the space~$\X_A$, which is the one-point compactification of the locally compact space
\begin{equation*}
\coprod_{n\in A}\Sph^{n+1}\amalg \omega.
\end{equation*}

We then identify the compact and Polish degree spectra of~$\X_A$ as follows.
\begin{theorem}\label{thm2}
For all~$A\subseteq\N$ and Turing degree~$\d$,
\begin{align*}
\text{$A$ is $\Sigma^0_2$ relative to~$\d$}&\iff\text{$\X_A$ has a $\d$-computable compact presentation,}\\
\text{$A$ is $\Sigma^0_3$ relative to~$\d$}&\iff\text{$\X_A$ has a $\d$-computable Polish presentation.}
\end{align*}
\end{theorem}
Before proving this result, let us give its main consequences.

\begin{corollary}\label{cor_one_jump}
For every Turing degree~$\d$, there exists a space~$\X_d$ whose compact degree spectrum is~$\{\x:\d\leq \x'\}$ and Polish degree spectrum is~$\{\x:\d\leq \x''\}$.
\end{corollary}
 \begin{proof}
 Given the Turing degree~$\d$ of a set~$D\subseteq\N$, let~$A=\{2n:n\in D\}\cup\{2n+1:n\notin D\}$ and~$\X_\d=\X_A$.
 \end{proof}

\begin{corollary}\label{cor_polish_nocompact}
There is a compact perfect Polish space which has a computable Polish presentation, but no computable compact presentation.
\end{corollary}

\begin{proof}
Let~$\X=\X_{\zero''}$ from Corollary \ref{cor_one_jump}. Its Polish degree spectrum is $\{\x:\zero''\leq\x''\}$ while its compact degree spectrum is $\{\x:\zero''\leq\x'\}$. The former contains~$\zero$ but not the latter.
\end{proof}

We can solve Question \ref{question:basic2} affirmatively.

\begin{corollary}
There is a perfect compact Polish space which is~$\zero'$-computably presentable but not computably presentable. Its Polish degree spectrum is the set of high$_2$ degrees.
\end{corollary}

\begin{proof}
Let~$X=\X_{\zero'''}$ from Corollary \ref{cor_one_jump}. Its Polish degree spectrum is~$\{\x:\zero'''\leq \x''\}$, i.e.~the high$_2$ degrees. It contains~$\zero'$ but not~$\zero$.
\end{proof}

One can improve the degree of non-computably presentable spaces.

\begin{corollary}
For any non-low$_2$ degree~$\d$, there exists a perfect Polish space which is $\d$-computably presentable, but not computably presentable.
\end{corollary}
\begin{proof}
The Polish degree spectrum of the space~$\X_{\d''}$ is~$\{\x:\d''\leq \x''\}$, so it contains~$\d$ but not~$\zero$ when~$\d$ is not low$_2$.
\end{proof}
In particular, one can take~$\d\leq \zero'$ which is low$_3$ but not low$_2$.

We now give the proof of Theorem \ref{thm2}. We need to measure the complexity of finding the dimension of the sphere. It turns out that the descriptive complexity of measuring its covering dimension is too large, as will be made precise later; the way we detect an~$n$-dimensional sphere is by finding a homology cycle of dimension~$n$. Although conceptually more complicated, it results in an algorithm of lower complexity.

\subsection{Learning the dimension of a sphere}
Given a compact presentation of space~$X$ which is homeomorphic to the~$d$-dimensional sphere~$\Sph^d$ for some~$d$, how difficult is it to find~$d$?
 
We are going to show that~$d$ can be computed in the limit, i.e.~one can compute a sequence~$(d_s)_{s\in\N}$ such that~$d=\lim_s d_s$.

The first idea would be to detect the covering dimension of~$X$, but it only enables to compute a sequence~$(d_s)_{s\in\N}$ such that~$d=\liminf_sd_s$, because the predicate~$\dim(X)\leq n$ is~$\Pi^0_2$ (Proposition \ref{prop_dim_pi02}). This strategy cannot work, because it can be proved that if~$X$ is homeomorphic to the $d$-dimensional cube~$[0,1]^d$, then one cannot compute~$d$ in the limit. Therefore, we need to use a way of distinguishing between spheres other than the dimension, and this is given by the homology of the spheres.

\begin{lemma}\label{lem_limit_sphere}
There is an algorithm that, given a compact presentation of a sphere~$\Sph^d$,~$d\geq 1$, computes~$d$ in the limit.
\end{lemma}
\begin{proof}
We show that relative to a compact presentation of~$X\cong \Sph^d$, the predicate~$X\cong \Sph^n$ is~$\Sigma^0_2$. It implies that it is~$\Delta^0_2$, because~$X\cong \Sph^n\iff\forall m\neq n,X\ncong\Sph^m$ is~$\Pi^0_2$, therefore~$\{d\}$ is~$\Delta^0_2$ and~$d$ can be computed in the limit.

The \v{C}ech homology groups of~$\Sph^d$ are
\begin{equation*}
\check{H}_n(\Sph^d;\zd)=\begin{cases}
\zd&\text{ if~$n=d$ or~$n=0$},\\
0&\text{otherwise.}
\end{cases}
\end{equation*}

Therefore, for~$n\geq 1$,~$X$ is homeomorphic to~$\Sph^n$ iff~$\check{H}_n(X;\zd)$ is non-trivial, which is~$\Sigma^0_2$ by Corollary \ref{cor_nontrivial_homology}.
\end{proof}
One could adapt the argument to include the case~$d=0$ in the previous lemma, but we will not need it.

Let us explain more concretely how the algorithm works. It searches for an approximation of~$X$ by a simplicial complex~$K_s$ that contains an~$n$-cycle~$c$ which is not a boundary, and which ``survives'' in each further approximation~$K_t$,~$t\geq s$. The right way of expressing that~$c$ survives in~$K_t$ is that there exists an~$n$-cycle~$d$ in~$K_t$ which ``refines'' $c$, in the sense that it is sent to~$c$ by the simplicial map~$f=f_{t-1}\circ\ldots \circ f_s:K_t\to K_s$, modulo a boundary: $c=f(d)+\partial e$, for some~$n$-chain~$e$ of~$K_s$.

Note that whether a cycle~$c$ of~$K_s$ survives at stage~$t$ is decidable, because it is expressed by  quantification over objects~$d,e$ in a finite set, and the maps~$f$ and~$\partial$ are computable.

At any stage~$t$, there are finitely many pairs~$(n,s)$ with~$s\leq t$ such that some~$n$-cycle in~$K_s$'s survives in~$K_t$. The algorithm outputs the number~$n$ for which the associated~$s$ is minimal, i.e.~$n$ is the dimension the cycle that has survived for the longest period of time.

\begin{remark}
Lemma \ref{lem_limit_sphere} can be alternatively proved by using the computability of \v{C}ech cohomology groups proved in \cite{LMN23} and \cite{DM23}. Indeed,~$\check{H}^n(\Sph^m)$ is non-trivial exactly when~$m=n$ (assuming~$m,n\geq 1$), and non-triviality of a~$\check{H}(\Sph^m)$ is~$\Sigma^0_2$.
\end{remark}

We can at last prove the main result of this section.
\begin{proof}[Proof of Theorem \ref{thm2}]
Assume that~$\d$ computes a compact presentation of~$\X_A$. Relative to~$\d$, we compute an enumeration~$(C_i)_{i\in\N}$ of the clopen subsets of~$\X_A$, by Proposition \ref{prop_clopen}. One has~$n\in A\iff\exists i,C_i$ is connected, is not a singleton and is homeomorphic to~$\Sph_{n+1}$. When~$C_i$ is connected and is not a singleton,~$C_i$ is a sphere so by Lemma \ref{lem_limit_sphere} its dimension is limit-computable relative to~$\d$. Therefore,~$A$ is~$\Sigma^0_2$ relative to~$\d$. If~$\d$ computes a Polish presentation of~$\X_A$, then~$\d'$ computes a compact presentation of~$\X_A$, so~$A$ is~$\Sigma^0_2$ relative to~$\d'$, hence~$\Sigma^0_3$ relative to~$\d$.

We now prove the other direction. If~$A$ is~$\Sigma^0_2$ relative to~$\d$, then there is a~$\d$-computable set~$B\subseteq \N^3$ such that~$n\in A\iff \exists i\forall j,(n,i,j)\in B$, and such that when~$n\in A$, there exists a unique~$i$ satisfying~$(n,i,j)\in B$ for all~$j$.

For each~$n$, let~$X^n_\infty$ be a canonical~$(n+1)$-sphere in~$\hil$ and~$X^n_j\subseteq X^n_\infty$ be an increasing sequence of finite sets at Hausdorff distance~$<2^{-j}$ from~$X^n_\infty$. We can take all these sets uniformly computably compact. Let
\begin{equation*}
Y^n_i=\begin{cases}
X^n_\infty&\text{if~$\forall j,(n,i,j)\in B$,}\\
X^n_j&\text{if~$j$ is minimal such that~$(n,i,j)\notin B$.}
\end{cases}
\end{equation*}

We apply Proposition \ref{prop_sequence_spaces}, implying that the sets~$Y^n_i$ are uniformly~$\d$-computably compact. Therefore, the set~$\alpha_0(\coprod_{n,i}Y^n_i\amalg\omega)$ is~$\d$-computably compact by Proposition \ref{prop_constructs}. This set is homeomorphic to~$\X_A$, because all the finite sets~$X^n_j$ can be seen as subsets of~$\omega$.

If~$A$ is~$\Sigma^0_3$ relative to~$\d$, then~$A$ can be presented in the same way for some~$\d$-c.e.~set~$B$, and we consider the same realization of~$\X_A$. Again by Proposition \ref{prop_sequence_spaces}, the sets~$Y^n_i$ are uniformly~$\d$-computably overt, and so is~$\alpha_0(\coprod_{n,i}Y^n_i\amalg\omega)$.
\end{proof}

\begin{remark}[Cycles vs dimension]
This result cannot be proved by simply using the dimension of the spheres, because detecting the dimension of a space is more complex in general. Whether a space has dimension~$n$ is not~$\Delta^0_2$; it is complete for the descriptive complexity class consisting of differences of two~$\Sigma^0_2$ sets. In particular, given a compact presentation of~$\I^n$, one cannot compute a sequence converging to~$n$. However, one can compute a sequence~$n_i$ such that~$n=\liminf_i n_i$, and this is optimal. It implies that if the space~$\X'_A:=\alpha_0(\coprod_{n\in A}\I^{n+1}\amalg\omega)$ has a computable compact presentation, then~$A$ is~$\Sigma^0_3$, and it is possible to find a~$\Sigma^0_3$-complete set~$A$ such that~$\X'_A$ has a computable compact presentation. However, it does not seem to be true that~$\X'_A$ has a computable compact presentation for any~$\Sigma^0_3$ set~$A$.
\end{remark}

 \begin{remark}
We do not know whether adding extra points outside the spheres is really needed. One can show that if~$Y_A$ is the one-point compactfication of~$\coprod_{n\in A}\Sph^{n+1}$, then~$Z$ computes a compact presentation of~$Y_A$ if and only if~$A$ is~$\Sigma^0_2$ relative to~$Z$. It implies that if~$Z$ computes a Polish presentation of~$Y_A$ then~$A$ is~$\Sigma^0_3$ relative to~$Z$, but we do not know whether the other implication holds.
 \end{remark}

\begin{remark}
One can modify the proof of Theorem \ref{thm2} to ensure that $\X_A$ is perfect, by using the one-point compactification of the following:
\begin{equation*}
\coprod_{n\in A}\mathbf{S}^{n+1}\amalg \coprod_{n\in\N}\I
\end{equation*}
where $\I$ is the line segment. As $\I$ has trivial homology groups, the algorithm that decodes~$A$ from~$X_A$ still works. In the construction of~$X_A$ from~$A$, one can replace the finite set of isolated points~$X^n_j$ with a finite disjoint union of line segments contained in~$\Sph^{n+1}$, and add countably many copies of $\I$.
\end{remark}

\subsection{Iterations of the Turing jump}
We now prove a result similar to Theorem \ref{thm2} for higher levels of the arithmetical hierarchy. This theorem also follows from the results independently obtained by Melnikov \cite{Melnikov21} and Lupini, Melinov, Nies  \cite{LMN23}, although with a different proof.

\begin{theorem}\label{thm_iteration}
Let~$k\geq 3$. To a set~$A\subseteq\N$ one can associate a compact Polish space~$\X_{A,k}$ such that for any Turing degree~$\d$,
\begin{align*}
\text{$A$ is $\Sigma^0_k$ relative to $\d$}&\iff\text{$\X_{A,k}$ has a~$\d$-computable compact presentation},\\
\text{$A$ is $\Sigma^0_{k+1}$ relative to $\d$}&\iff\text{$\X_{A,k}$ has a~$\d$-computable Polish presentation}.
\end{align*}
\end{theorem}
\begin{corollary}\label{cor_iteration}
For every Turing degree~$\d$ and every~$k\geq 2$, there exists a space~$\X_{\d,k}$ whose compact degree spectrum is~$\{\x:\d\leq_T \x^{(k)}\}$ and Polish degree spectrum is~$\{\x:\d\leq_T \x^{(k+1)}\}$.
\end{corollary}

\begin{proof}
 Given the Turing degree~$\d$ of a set~$D\subseteq\N$, let~$A=\{2n:n\in D\}\cup\{2n+1:n\notin D\}$ and~$\X_{\d,k}=\X_{A,k+1}$, and note that for~$j\in\{k,k+1\}$,~$A$ is~$\Sigma^0_{j+1}$ relative to~$\x$ iff~$\d\leq \x^{(j)}$.
\end{proof}

We recall that a degree~$\x$ is high$_n$ iff~$\zero^{(n+1)}\leq \x^{(n)}$. The class of high$_n$ degrees is the compact degree spectrum of the computable Polish space~$\X_{\zero^{(n+1)},n}$ and is the Polish degree spectrum of the space~$\X_{\zero^{(n+1)},n-1}$.

We now proceed with the proof of Theorem \ref{thm_iteration}. We restrict our attention to compact Polish spaces that have countably many connected components. For such spaces, the dimension is the supremum of the dimensions of the connected components (sum theorem 1.5.3 in \cite{EngDim78}).


\paragraph{Widgets.} For~$k\in\N$, we inductively define spaces~$S_k$ and~$P_k$.

For a space~$X$, let~$\alpha_0(X)$ be the one-point compactification of the disjoint union of countably many copies of~$X$. It will be thought as a pointed space, with the point at infinity as basepoint. For~$k\geq 1$ and a compact space~$X$, let~$\alpha_k(X)$ be the wedge sum of the pointed space~$(\I^k,(0,\ldots,0))$ with~$\alpha_0(X)$.

 Let~$S_0=\emptyset$ and~$P_0=\{0\}$. For~$k\geq 1$, let
\begin{align*}
S_k&=\alpha_k(S_{k-1}\amalg P_{k-1}),\\
P_k&=\alpha_k(S_{k-1}).
\end{align*}

In particular,
\begin{align*}
S_1&=\I\vee (\omega+1),\\
P_1&=\I.
\end{align*}

For~$k\geq 1$, these spaces have several properties:
\begin{itemize}
\item $\dim(S_k)=\dim(P_k)=k$,
\item If~$C$ is a clopen subset of~$S_k$ (respectively~$P_k$) of dimension~$k$, then~$C$ contains a clopen set that is homeomorphic to~$S_k$ (respectively~$P_k$).
\end{itemize}

\paragraph{A family of topological invariants.}
For~$k\geq 1$, we define a topological invariant that distinguishes~$S_k$ from~$P_k$ and whose complexity in the Vietoris topology is~$\Sigma^0_k$, in the sense that the collection of compact subsets of~$\hil$ satisfying this invariant is~$\Sigma^0_k$ in the hyperspace~$\K(\hil)$ endowed with the Vietoris topology. More generally, we define a~$\Sigma^0_k$ invariant~$\J_{a,k}$ that distinguishes~$\I^a\times S_k$ from~$\I^a\times P_k$.

\begin{definition}For~$a\in \N$ and~$k\geq 1$, we inductively define a topological invariant~$\J_{a,k}$ as follows:
\begin{enumerate}
\item $X\in\J_{a,1}\iff X$ is disconnected,
\item For~$k\geq 2$,~$X\in \J_{a,k}\iff X$ contains a clopen set of dimension~$\geq a+k-1$ that does not satisfy~$\J_{a,k-1}$.
\end{enumerate}
\end{definition}

\begin{remark}
$\J_{a,k}$ enjoys some form of heredity:~$X$ satisfies~$\J_{a,k}$ iff~$X$ contains a clopen set satisfying~$\J_{a,k}$.
\end{remark}

\begin{lemma}
Let~$k\geq 1$. $\J_{a,k}$ is satisfied by~$\I^a\times S_k$ but not by~$\I^a\times P_k$.
\end{lemma}
\begin{proof}
For~$k=1$, it is straightforward:~$\I^a\times S_1=\I^a\times (\I\vee (\omega+1))$ is disconnected, but~$\I^a\times P_1=\I^{a+1}$ is connected.

Let~$k\geq 2$. We assume the result for~$k-1$, and prove it for~$k$.~$\I^a\times S_k$ contains a clopen set~$C$ which is a copy of~$\I^a\times P_{k-1}$. One has~$\dim(C)=a+k-1$ and by induction hypothesis, this set does not satisfy~$\J_{a,k-1}$.

Let~$C$ be a clopen subset of~$I^a\times P_k$ of dimension~$\geq a+k-1$. We show that~$C$ satisfies~$\J_{a,k-1}$. We claim that~$C$ contains a clopen set~$D$ which is a copy of~$\I^a\times S_{k-1}$, which implies that~$C$ satisfies~$\J_{a,k-1}$ by heredity.

If~$\dim(C)=a+k$ then~$C$ contains a clopen set which is a copy of~$\I^a\times P_k$, which in turns contains a clopen copy of~$\I^a\times S_{k-1}$. If~$\dim(C)=a+k-1$ then~$C$ is disjoint from~$\I^a\times \I^k$, so it is a disjoint union of clopen subsets of~$\I^a\times S_{k-1}$. One of them, call it~$E$, must have dimension~$a+k-1$. Being a clopen subset of~$\I^a\times S_{k-1}$ of dimension~$a+k-1$,~$E$ contains a clopen set~$D$ which is a copy of~$\I^a\times S_{k-1}$. $D$ is a clopen subset of~$C$ as wanted.

Therefore, we have proved that~$C$ satisfies~$\J_{a,k-1}$ for every clopen subset of~$\I^a\times P_k$ of dimension~$\geq a-k-1$, so~$\I^a\times P_k$ does not satisfy~$\J_{a,k}$.
\end{proof}

\begin{lemma}
For~$k\geq 1$,~$\J_{a,k}$ has complexity~$\Sigma^0_k$ in the Vietoris topology.
\end{lemma}
\begin{proof}
We prove the result by induction. For~$k=1$, being disconnected is~$\Sigma^0_1$.

Let~$k\geq 2$ and assume that~$\J_{a,k-1}$ is~$\Sigma^0_{k-1}$. Having dimension~$\geq a+k-1$ is~$\Sigma^0_2$ (Proposition \ref{prop_dim_pi02}) hence~$\Sigma^0_k$, the complement of~$\J_{a,k-1}$ is~$\Pi^0_{k-1}$ by induction hypothesis. Therefore, containing a clopen set of dimension~$\geq a+k-1$ that does not satisfy~$\J_{a,k-1}$ is~$\Sigma^0_k$.
\end{proof}

\paragraph{Encoding a set in a space.}
Let~$k\geq 3$ be fixed. We define~$a_n=(k+1)n$ for all~$n\in\N$.

Let~$A\subseteq \N$. We define the space~$\X_{A,k}$ as the one-point compactification of the disjoint union of the spaces~$\I^{a_n}\times S_k$ for~$n\in A$ and~$\I^{a_n}\times P_k$ for~$n\notin A$.

The role of~$a_n$ is to use the dimension to distinguish between the difference pieces. In particular, if a clopen subset of~$\X_{A,k}$ has dimension~$< a_n$, then~$C$ does not intersect~$\I^{a_n}\times S_k$ or~$\I^{a_n}\times P_k$.
\begin{lemma}
One has~$n\in A$ iff~$\X_{A,k}$ contains a clopen set of dimension~$\leq a_n+k$ satisfying~$\J_{a_n,k}$.
\end{lemma}
\begin{proof}
Let~$n\in\N$ and~$a=a_n$.

First assume that~$n\in A$. $\X_{A,k}$ contains~$\I^{a}\times S_k$ as a clopen set, which has dimension~$a+k$ and satisfies~$\J_{a,k}$.

Now let~$n\notin A$ and assume for a contradiction that~$\X_{A,k}$ contains a clopen set~$C$ of dimension~$\leq a+k$ satisfying~$\J_{a,k}$. First,~$C$ is contained in the union of the sets~$\I^{a_m}\times S_k$ or~$\I^{a_m}\times P_k$ for~$m\leq n$. As~$C$ satisfies~$\J_{a,k}$,~$C$ contains a clopen set~$D$ of dimension~$\geq a+k-1$ that does not satisfy~$\J_{a,k-1}$. Therefore, the intersection~$E=D\cap (\I^a\times P_k)$ does not satisfy~$\J_{a,k-1}$ by heredity. One has~$\dim(E)\geq a+k-1$, because the part of~$D$ which is outside~$\I^a\times P_k$ has dimension~$<a$. As a result,~$\I^a\times P_k$ satisfies~$\I^a\times P_k$, which is a contradiction.
\end{proof}

Note that the property expressed in the lemma is~$\Sigma^0_k$ in the Vietoris topology, when~$k\geq 3$ (however it is not~$\Sigma^0_2$ when~$k=2$, so Theorem \ref{thm2} cannot be proved in this way). Therefore, for any Turing degree~$\d$ computing a compact presentation of~$\X_{A,k}$,~$A$ is~$\Sigma^0_k$ relative to~$\d$; hence for any Turing degree~$\d$ computing a Polish presentation of~$\X_{A,k}$,~$\d'$ computes a compact presentation of~$\X_{A,k}$ (Lemma \ref{lem:Polish-to-compact}) and therefore~$A$ is~$\Sigma^0_{k+1}$ relative to~$\d$.

We now show that if~$A$ is~$\Sigma^0_k$ relative to~$\d$, then~$\d$ computes a compact presentation of~$\X_{A,k}$, and if~$A$ is~$\Sigma^0_{k+1}$ relative to~$\d$, then~$\d$ computes a Polish presentation of~$\X_{A,k}$. Thanks to Proposition \ref{prop_sequence_spaces}, the two results can be proved with a single argument. We prove it by induction on~$k$, the next result is the case~$k=1$.
\begin{lemma}\label{lem_S1P1}
For a set~$A$, let
\begin{equation*}
Y_n=\begin{cases}
S_1&\text{if }n\in A,\\
P_1&\text{if }n\notin A.
\end{cases}
\end{equation*}
If~$A$ is~$\Sigma^0_1$ (resp.~$\Sigma^0_2$) relative to~$\d$, then there exist copies of~$Y_n$ that are uniformly~$\d$-computably compact (resp.~overt).
\end{lemma}

\begin{proof}
For~$k\in\N$, let
\begin{align*}
X_k&=[2^{-k},1]\cup\{2^{-k}-2^{-n}:n\geq k\},\\
X_\infty&=[0,1].
\end{align*}
Observe that~$X_k\cong S_1$,~$X_\infty\cong P_1$,~$X_k\subseteq X_{k+1}\subseteq\ldots\subseteq X_\infty$ and~$d_H(X_k,X_\infty)<2^{-k-1}$.

Let~$A$ be given by a predicate~$n\in A\iff \exists k,(n,k)\in B$. To~$n$ we associate~$Y_n=X_k$ if~$k$ is minimal such that~$(n,k)\in B$, and~$Y_n=X_\infty$ if there is no such~$k$, i.e.~$n\notin A$. 

We apply Proposition \ref{prop_sequence_spaces}. If~$B$ is computable relative to~$\d$, then~$Y_n$ is unifomly~$\d$-computably compact. If~$B$ is~$\Pi^0_1$ relative to~$\d$, then~$Y_n$ is uniformly~$\d$-computably overt.
\end{proof}

\begin{lemma}
Let~$k\geq 1$. For a set~$A$, let
\begin{equation*}
Y_n=\begin{cases}
S_k&\text{if }n\in A,\\
P_k&\text{if }n\notin A.
\end{cases}
\end{equation*}
If~$A$ is~$\Sigma^0_k$ (resp.~$\Sigma^0_{k+1}$) relative to~$\d$,  then there exist copies of~$Y_n$ that are uniformly~$\d$-computably compact (resp.~overt).
\end{lemma}
\begin{proof}
The case~$k=1$ is Lemma \ref{lem_S1P1}. Let~$k\geq 2$ and assume the result for~$k-1$.

Let~$A$ be given by a predicate~$n\in A\iff \exists i,(n,i)\in B$ where~$B$ is~$\Pi^0_{k-1}$ relative to~$\d$. We can assume that if~$\exists i,(n,i)\in B$, then there exist infinitely many~$i$'s such that~$(n,i)\in B$ and infinitely many~$i$'s such that~$(n,i)\notin B$. Indeed, define~$B'$ as follows:~$(n,2i)\in B'$ if there exists~$j\leq i$ such that~$(n,j)\in B$, and~$(n,2i+1)\notin B'$. $B'$ is~$\Pi^0_{k-1}$ and~$n\in A\iff\exists i,(n,i)\in B'$.

Let~$X_{n,i}=S_{k-1}$ if~$(n,i)\notin B$ and~$X_{n,i}=P_{k-1}$ if~$(n,i)\in B$. By induction hypothesis, given~$(n,i)$ one can produce a copy of~$X_{n,i}$ which is uniformly~$\d$-computable compact. Build the space~$X_n$ which is the wedge sum of~$I^k$ and the one-point compactification of~$\coprod_{n,i}X_{n,i}$. If~$n\in A$ then~$X_n=S_k$; if~$n\notin A$ then~$X_n=P_k$.

If~$B$ is~$\Pi^0_k$ relative to~$\d$, then the very same argument gives uniformly~$\d$-computably overt copies of~$Y_n$ by induction.
\end{proof}

Finally, a presentation of~$\X_{A,k}$ from~$A$ is then achieved by building products of~$\I^{a_n}$ with copies of~$S_k$ or~$P_k$ as above, using Proposition \ref{prop_constructs}. The proof of Theorem \ref{thm_iteration} is complete.

\section{Cone-avoidance}\label{sec_cone}

In this section, we solve Question \ref{question:basic3} by showing the following:

\begin{theorem}\label{thm_cone}
Let~$A\subseteq\N$ be a non-c.e.~set. Every compact Polish space has a Polish presentation that does not enumerate~$A$.
\end{theorem}

In particular, 
\begin{corollary}
The degree spectrum of a compact Polish space cannot be the upper cone~$\{\x:\d\leq \x\}$ for any non-computable degree~$\d$.
\end{corollary}

Actually the proof also shows that if for each~$i$ we choose a non-c.e.~set~$A_i$, then every compact Polish space has a presentation that does not enumerate any~$A_i$. It implies that the degree spectrum of a compact Polish space cannot be a countable union of non-trivial upper cones~$\bigcup_{i\in\omega}\{\x:\d_i\leq\x\}$.

In order to prove the result, we need ideas from computability theory and ideas from topology.

\paragraph{Overtness argument.}
Overtness captures a familiar argument in computably theory, which we describe now.

We will apply the technique to the space~$X=\V(Q)$, however it is easier to state the result for an abstract space~$X$.

Let~$X$ be a countably-based space with a fixed indexed basis~$(B_i)_{i\in\N}$ that is closed under finite intersections. We say that~$A\subseteq\N$ is reducible to~$x\in X$, written~$A\leq_e x$, if~$A$ is enumeration reducible to~$N_x=\{i\in\N:x\in B_i\}$. We write~$M(x)=A$ if~$M$ enumerates~$A$ from any enumeration of~$N_x$. We say that a set~$S\subseteq X$ is computably overt if the set~$\{i\in\N:S\cap B_i\neq\emptyset\}$ is c.e.

For a Turing machine~$M$ and a set~$A\subseteq\N$, we say that~$M$ fails to enumerate~$A$ from~$x$ if~$M$ outputs some~$n\notin A$ on some enumeration of~$N_x$. We denote by~$F_{M,A}$ the set of~$x$'s on which~$M$ fails to enumerate~$A$. That set is open, because when~$M$ outputs some~$n\notin A$, it has only read a finite part of~$N_x$, which can be extended to an enumeration of~$N_y$ for any~$y$ in some neighborhood of~$x$, so each such~$y$ also belongs to~$F_{M,A}$.

\begin{lemma}\label{lem_overt}
Let~$X$ be a countably-based space. Let~$x\in X$,~$A\subseteq\N$ be a non-c.e.~set and~$M$ a Turing machine. If~$M(x)=A$, then~$F_{M,A}$ intersects every computably overt subset of~$X$ containing~$x$.
\end{lemma}
\begin{proof}
Assume that~$M(x)=A$ and let~$V$ be a computably overt set containing~$x$. If~$F_{M,A}$ does not intersect~$V$, then we describe an effective procedure that enumerates~$A$, contradicting the assumption that~$A$ is not c.e. The procedure is as follows: enumerate all the prefixes of names of elements of~$V$ (which is possible because~$V$ is computably overt), simulate~$M$ on them, and collect all the outputs. As~$F_{M,A}$ does not intersect~$V$, all the outputs are correct, i.e.~belong to~$A$, and every element of~$A$ appears because~$M$ enumerates~$A$ on each name of~$x\in V$. As a result, this procedure enumerates~$A$, which is then c.e. The contradiction implies that~$F_{M,A}$ intersects~$V$.
\end{proof}

\paragraph{Perturbations.}
We now come to the topological ingredient of the proof, based on the notion of~$\epsilon$-perturbation. The proof is a Baire category argument: if~$A\subseteq\N$ is not c.e., then one can perturb any compact set~$K\subseteq\hil$ so that its perturbed copy~$K'$ does not enumerate~$A$, when seeing~$K'$ as an element of the space~$\V(\hil)$.

\begin{definition}
An~$\epsilon$\textbf{-perturbation} is a one-to-one continuous function~$f:\hil\to\hil$ such that~$d(f(x),x)<\epsilon$ for all~$x\in\hil$.
\end{definition}

\begin{lemma}\label{lem_finite_sets}
Let~$S=\{s_0,\ldots,s_n\}$ and~$T=\{t_0,\ldots,t_n\}$ be finite subsets of~$\hil$ such that~$d(s_i,t_i)<\epsilon$ for~$i=0,\ldots ,n$. There exists an~$\epsilon$-perturbation~$f$ such that~$f(s_i)=t_i$ for~$i=0,\ldots,n$.
\end{lemma}
\begin{proof}
It is a direct application of the homeomorphism extension theorem from \cite{vM01}. It relies on the notion of~$Z$-set, which we recall for completeness; however, it is not important to understand this notion, we will only use the fact that finite subsets of~$\hil$ are~$Z$-sets.

The following definition and results are taken from \cite{vM01}. A closed set~$A\subseteq\hil$ is a~$Z$-set if for every continuous function~$f:\hil\to\hil$ and every~$\epsilon>0$, there exists a continuous function~$g:\hil\to\hil$ such that~$d(f(x),g(x))<\epsilon$ for all~$x\in\hil$, and~$g(Q)\cap A=\emptyset$ (see \S 5.1 in \cite{vM01}). By Remark 5.1.4 in \cite{vM01}, every singleton is a~$Z$-set. Lemma 5.1.2 (3) in \cite{vM01} states that any finite union of~$Z$-sets is a~$Z$-set, therefore every finite set is a~$Z$-set. 

The homeomorphism extension theorem (Theorem 5.3.7 in \cite{vM01}) states that if~$S,T\subseteq \hil$ are~$Z$-sets and~$f:S\to T$ is a homemorphism satisfying~$d(f(x),x)<\epsilon$ for all~$x\in S$, then~$f$ can be extended to a homeomorphism~$\overline{f}:\hil\to \hil$ satisfying~$d(\overline{f}(x),x)<\epsilon$ for all~$x\in\hil$. Therefore, the statement is just an application of the homeomorphism extension theorem to finite sets.
\end{proof}

We remind the reader that~$\V(\hil)$ is a countably-based topological space endowed with the lower Vietoris topology. In the next statement, the notions of computable overtness and closure are meant in that topology.

\begin{lemma}\label{lem2}
Let~$K\subseteq \hil$ be a compact set and~$\epsilon>0$. There exists a computably overt set~$\A\subseteq\V(\hil)$ with~$K\in \A$, such that~$\A$ is contained in the closure of the set of~$\epsilon$-deformations of~$K$.
\end{lemma}
\begin{proof}
Let~$\F\subseteq\V(\hil)$ be the family of finite sets of rational points, which can be computably indexed in an obvious way. We are going to define~$\A$ in such a way that:
\begin{enumerate}
\item $\F\cap\A$ is dense in~$\A$ (w.r.t.~the lower Vietoris topology),
\item $\F\cap\A$ is computably enumerable,
\item Every element of~$\F\cap\A$ is contained in some~$\epsilon$-deformation of~$K$.
\end{enumerate}
The first two conditions imply that~$\A$ is computably overt, because it is the closure of a computable sequence listing the elements of~$\F\cap \A$.

The first and third conditions imply that~$\A$ is contained in the closure of the set of~$\epsilon$-deformations of~$K$:~$\A$ is the closure of~$\F\cap\A$, and each element of~$\F\cap\A$ is a subset of an~$\epsilon$-deformation~$K'$ of~$K$, so belongs to the closure of~$\{K'\}$ (the lower Vietoris closed open sets are upwards closed, equivalently the lower Vietoris closed sets are downwards closed).

We now define~$\A$ satisfying these conditions. If~$K$ was perfect then we could just take some small rational ball in the Hausdorff metric containing~$K$. However we need a bit more work in the general case.

We first show that there exist open sets~$U_0,\ldots,U_n\subseteq\hil$ that cover~$K$, such that for every~$x\in U_0$,~$K\cap B(x,\epsilon)$ is infinite and~$K\cap U_i$ is a singleton for each~$i\geq 1$. Let~$K_N$ be the set of non-isolated points of~$K$ and let~$U_0=\{x\in Q:d(x,K_N)<\epsilon\}$. The set~$K\setminus U_0$ is finite, because it is compact and all its points are isolated. Therefore, there exist basic balls~$U_1,\ldots,U_n$ isolating the points of~$K\setminus U_0$. We can make sure that all the open sets~$U_i$ are pairwise disjoint. By compactness of~$K$, we can assume that~$U_0$ is a finite union of basic balls. We can now define our computably overt subset of~$\V(\hil)$: let
\begin{equation*}
\A=\left\{C\in\V(\hil):C\subseteq \bigcup_{0\leq i\leq n} U_i\text{ and }|C\cap U_i|=1\text{ for all }i\geq 1\right\}.
\end{equation*}
We check condition 1. Let~$C\in\A$. $C$ is a limit in the Vietoris topology of finite sets~$C_k$ of rational points (note that here we use the Vietoris rather than the lower Vietoris topology). As~$C\in \A$, for sufficiently large~$k$,~$C_k$ is contained in~$\bigcup_i U_i$ and intersects~$U_1,\ldots,U_n$. For~$i\geq 1$, if~$C_k\cap U_i$ contains more than one point, then keep only one of them. Let~$C'_k\subseteq C_k$ be obtained this way. One has~$C'_k\in\F\cap \A$ by construction, and~$C'_k$ converges to~$C$ in the Vietoris (hence lower Vietoris) topology. 

Condition 2.~is easily checked: the conditions defining~$C\in \A$ are c.e., when~$C$ is a finite set of rational points.




We now prove condition 3. We show that for every finite set~$T\in \A$, there exists an~$\epsilon$-deformation of~$K$ containing~$T$.

Let~$T=\{t_0,\ldots,t_n\}$ belong to~$\A$. We build a finite set~$S=\{s_0,\ldots,s_n\}\subseteq K$ with~$d(s_i,t_i)<\epsilon$. For each~$i$, if~$t_i\in U_0$ then the intersection of~$B(t_i,\epsilon)$ with~$K$ is infinite, so we can choose a point~$s_i$ in the intersection, so that~$s_i\neq s_j$ for~$i\neq j$. If~$t_i\in U_k$ with~$k\geq 1$, then we define~$s_i$ as the unique point of~$K$ in~$U_k$. The points~$s_i$ are pairwise distinct, because if~$i\neq j$ then~$s_i$ and~$s_j$ cannot both belong to a common~$U_k$,~$k\geq 1$, as~$K\in \A$.

One has~$d(s_i,t_i)<\epsilon$ for each~$i$, so we can apply Lemma \ref{lem_finite_sets} to obtain an~$\epsilon$-perturbation~$f$ mapping each~$s_i$ to~$t_i$. One has~$T=f(S)\subseteq f(K)$ so the proof is complete.
\end{proof}

We now have all the ingredients needed to prove the result.
\begin{proof}[Proof of Theorem \ref{thm_cone}]
Let~$X$ be a compact Polish space. We prove that some copy of~$X$ in~$\V(\hil)$ does not enumerate~$A$. It implies the result, because any name of a copy of~$X$ in~$\V(\hil)$ computes a presentation of~$X$.

The space~$\inj(\hil)=\{\phi:\hil\to\hil\text{ continuous one-to-one}\}$, with the topology induced by the metric~$\rho(\phi,\psi)=\sup_x d(\phi(x),\psi(x))$ is Polish (Corollary 1.3.11 in \cite{vM01}). For any compact set~$K\subseteq\hil$ and any non-c.e.~set~$A\subseteq\N$, we prove that the set~$\{\phi\in\inj(\hil):A\leq_e \phi(K)\}$ is meager in~$\inj(\hil)$, which implies the existence of a copy~$K'$ of~$K$ such that~$A\nleq_e K'$ ($K'$ is seen as an element of the space~$\V(\hil)$). It is done by showing that for each Turing machine~$M$, the set~$\{\phi\in\inj(\hil):M(\phi(K))=A\}$ is nowhere dense in~$\inj(\hil)$.

Let~$\phi\in\inj(\hil)$ be such that~$M(\phi(K))=A$. For~$\epsilon>0$, we prove that there exists~$\psi\in\inj(\hil)$ such that~$\rho(\phi,\psi)<\epsilon$ and such that~$\psi(K)\in F_{M,A}$ (the set~$F_{M,A}$ was defined before Lemma \ref{lem_overt}). It implies the result, because for every~$\psi'$ sufficiently close to~$\psi$, one also has~$\psi'\in F_{M,A}$ as~$F_{M,A}$ is open.

Lemma \ref{lem2} provides a computably overt set~$\A$ containing~$\phi(K)$, in which the set of~$\epsilon$-deformations of~$\phi(K)$ is dense. By Lemma \ref{lem_overt},~$F_{M,A}$ intersects~$\A$. As~$F_{M,A}$ is open, there exists an~$\epsilon$-deformation of~$\phi(K)$ in~$F_{M,A}$. Let~$f$ be the corresponding~$\epsilon$-perturbation, and~$\psi=f\circ \phi$. One has~$\rho(\psi,\phi)<\epsilon$ and~$\psi(K)\in F_{M,A}$.
\end{proof}

\section{Comparing compact and Polish presentations}\label{sec:jump}

Let~$X$ be a compact Polish space. In the proofs we frequently used the fact that the jump of any Polish presentation of~$X$ computes a compact presentation of~$X$, which is stated precisely in Lemma \ref{lem:Polish-to-compact}. Here we investigate whether it can compute more. Of course, it always computes~$\zero'$, and we show that if~$X$ is perfect, then it does not compute more in general: every compact presentation of~$X$, paired with~$\zero'$, computes the jump of a Polish presentation of~$X$. However, we will see after that it is no more true for non-perfect spaces and we give a counter-example.

\begin{theorem}\label{thm_jump}
Let~$X$ be a compact perfect Polish space. If~$X$ has a~$\d$-computable compact presentation, then there exists~$\e$ such that~$X$ has an~$\e$-computable Polish presentation and~$\e'\leq \d\oplus \zero'$.
\end{theorem}

\begin{corollary}\label{cor_degree_perfect}
Let~$X$ be a compact perfect Polish space. One has
\begin{equation*}
\{\e':\e\in\Spec(X)\}=\{\d\oplus \zero':\d\in\Spec_c(X)\}.
\end{equation*}
\end{corollary}
\begin{proof}
If~$\e\in\Spec(X)$, then let~$\d=\e'$. One has~$\d\in \Spec_c(X)$ by Lemma \ref{lem:Polish-to-compact}, and~$\e'=\d\oplus\zero'$.

If~$\d\in\Spec_c(X)$, then there exists~$\e_0\in\Spec(X)$ such that~$\e_0'\leq \d\oplus \zero'$ by Theorem \ref{thm_jump}. Friedberg's jump inversion theorem relative to~$\e_0$ states that if~$\mathbf{a}\geq \e_0'$, then there exists~$\mathbf{f}$ such that~$\mathbf{a}=(\e_0\oplus \mathbf{f})'$. We apply it to~$\mathbf{a}=\d\oplus \zero'$, and let~$\e=\e_0\oplus \mathbf{f}$. One has~$\d \oplus\zero'=\e'$ and~$\e\in\Spec(X)$ as~$\e\geq \e_0\in\Spec(X)$.
\end{proof}

\paragraph{Reformulation.}
Again, we use an overtness argument to reformulate the problem.


As we have already seen, the degrees of Polish presentations of~$X$ coincide with the degrees of copies of~$X$ as elements of~$\V(\hil)$. In the same way, we show that the jumps of these degrees are exactly the jumps of the copies of~$X$ in~$\V(\hil)$.

Again, we abstract away from~$\V(\hil)$ to make the results easier to read. Let~$S$ be a countably-based space with a fixed index basis~$(B_i)_{i\in\N}$ that is closed under finite intersections. Let~$(U_i)_{i\in\N}$ be the canonical enumeration of the effective open subsets of~$S$ defined by~$U_i=\bigcup_{j\in W_i}B_j$, where~$W_i$ is the~$i$th c.e.~subset of~$\N$. We will apply the next result to~$S=\V(\hil)$.

\begin{definition}
The jump of~$x\in S$ is the set~$J(x)=\{i\in\N:x\in U_i\}$.
\end{definition}
\begin{lemma}
If~$x$ is a point of an effective countably-based space~$S$, the Turing degree of~$J(x)$ is the least element of~$\{\d':\d\text{ computes }x\}$.
\end{lemma}
\begin{proof}
Let~$\delta_S:\subseteq\Baire\to S$ be the standard representation of~$S$, mapping~$p$ to~$x$ if~$\{i:\exists n,p(n)=i+1\}=\{i:x\in B_i\}$, in which case we say that~$p$ is a name of~$x$. The function~$\delta_S$ is computable and effectively open: the image of an effective open set is an effective open set, uniformly. Observe that~$\d$ computes~$x$ if and only if~$\d$ computes some name of~$x$.

As~$\delta_S$ is computable, the preimages of effective open sets are effectively open, so if~$p$ is a name of~$x$ then~$p'$ computes~$J(x)$.

Conversely, we show that~$J(x)$ computes~$p'$ for some name~$p$ of~$x$. Let~$(U_n)_{n\in\N}$ be the canonical enumeration of the effective open subsets of~$\Baire$. The set~$\delta_S^{-1}(x)$ is~$\Pi^0_2$ relative to~$J(x)$: indeed, one has~$p\in\delta_S^{-1}(x)$ iff~$\forall i,[x\in B_i\iff \exists n,p(n)=i+1]$, which is a~$\Pi^0_2$ formula relative to~$J(x)$. Let then~$V_n$ be uniformly~$J(x)$-effective open sets such that~$\delta_S^{-1}(x)=\bigcap_nV_n$. Moreover, if~$W\subseteq\Baire$ is an effective open set, given by a index, then we can decide using~$J(x)$ whether~$W$ intersects~$\delta_S^{-1}(x)$, because it is equivalent to~$x\in \delta_S(W)$, which is an effective open set for which we have an index.

At stage~$s$, we have produced a finite prefix~$p_s$ of~$p$, such that~$[p_s]$ is contained in~$V_{s-1}$ and intersects~$\delta^{-1}_S(x)$. Using~$J(x)$, we can decide whether~$[p_s]\cap U_s$ intersects~$\delta_S^{-1}(x)$. If it does, then we extend~$p_s$ to~$p_{s+1}$ so that~$[p_{s+1}]$ is contained in~$U_s\cap V_s$ and intersects~$\delta_S^{-1}(x)$, and we declare that~$p\in U_s$. If it does not, then we simply extend~$p_s$ to~$p_{s+1}$ so that~$[p_{s+1}]$ is contained in~$V_s$ and intersects~$\delta_S^{-1}(x)$, and we declare that~$p\notin U_s$.

In the limit, we obtain some~$p\in\bigcap_s V_s$ so~$p$ is a name of~$x$. For each~$s$, we have decided along the construction whether~$p\in U_s$, so we have computed~$p'$.
\end{proof}

\begin{proof}[Proof of Theorem \ref{thm_jump}]
Let~$X\in\K(\hil)$ be~$\d$-computably compact and perfect. Using~$\d$ and~$\zero'$ as oracles, we progressively compute a copy~$K$ of~$X$, together with its jump~$J(K)$, seeing~$K$ as a point of~$\V(\hil)$. Let~$(\U_n)_{n\in\N}$ be an effective enumeration of the effective open subsets of~$\V(\hil)$. For each~$n$, we need to decide as long as we build~$K$, whether~$K\in \U_n$.

We start from some~$\epsilon_0>0$ and some basic~$\epsilon_0/2$-ball~$\B_0$ in the Hausdorff metric, containing~$X$. $\B_0$ is a computably overt set in~$\V(\hil)$, so we can decide using~$\zero'$ whether it intersects~$\U_0$. There are two cases:

\begin{enumerate}
\item Assume that~$\B_0$ intersects~$\U_0$. As in Lemma \ref{lem2} there exists an~$\epsilon_0$-perturbation~$f_0$ mapping~$X$ to~$X_0\in \U_0$ (as~$X$ is perfect, the computably overt set given by Lemma \ref{lem2} can be replaced by the~$\epsilon_0/2$-ball~$\B_0$).

We now show that one can compute such an~$f_0$. The space~$\P_{\epsilon_0}$ of~$\epsilon_0$-perturbations is a computable Polish space, the function~$\Phi:\P_{\epsilon_0}\to \V(\hil)$ mapping~$f$ to~$f(X)$ is~$\d$-computable, so~$\Phi^{-1}(\U_0)$ is a~$\d$-effective open set, therefore one can~$\d$-computably find some~$f_0$ there.

We now pick a ball~$\B_1$ of some radius~$\epsilon_1$ (to be chosen, see below) around~$X_0$, whose closure is contained in~$\B_0\cap \U_0$ and in which we are going to stay forever, so that in the limit, the copy~$K$ of~$X$ belongs to~$\U_0$. We declare that~$K\in \U_0$.

\item Now assume that~$\B_0$ does not intersect~$\U_0$. We do nothing and proceed by picking a ball~$\B_1$ of some radius~$\epsilon_1$ around~$X$ whose closure is contained in~$\B_0$. In the sequel, we stay forever in this ball so in the limit,~$K\in \B_0$ hence~$K$ does not belong to~$\U_0$. We declare that~$K\notin \U_0$.
\end{enumerate}

In both cases, we have decided whether the set~$K$ belongs to~$\U_0$. We now iterate this process with~$\U_1,\U_2$, etc., choosing~$\epsilon_1,\epsilon_2,\ldots$ smaller and smaller so that the composition of the~$\epsilon_i$-perturbations converges to an injective function (the Inductive Convergence Criterion \cite[Theorem 1.6.2]{vM01} tells us that we can always choose the next~$\epsilon_i$ sufficiently small to ensure that the limit is a injective, and moreover~$\epsilon_i$ can be chosen in a computable way), and taking the closure of~$\B_{n+1}$ contained in~$\B_n$. In the limit, we have built a copy~$K$ of~$X$ and computed its jump, using~$\d$ and~$\zero'$ as oracles.
\end{proof}

The argument is uniform: there is a uniform procedure relative to~$\zero'$ that, given a compact presentation of any compact perfect Polish space~$\X$, computes the jump of some Polish presentation of~$\X$. We now observe that there cannot exist a uniform argument including non-perfect Polish spaces. Indeed, whether~$X$ is not perfect is~$\Sigma^0_2$ for a Polish presentation, so it is~$\Sigma^0_1$ in its jump. If there was a uniform argument relative~$\zero'$, then being non-perfect would be~$\Sigma^0_1(\zero')$ for compact presentations, in particular it would be open in~$\K(Q)$ endowed with the Hausdorff metric. However, the set of non-perfect compact sets is not open, as witnessed for instance by a sequence of segments shrinking to a singleton.

We actually show that Theorem \ref{thm_jump} simply does not extend to non-perfect Polish spaces. Observe that a consequence of Theorem \ref{thm_jump} is that if~$X$ is a perfect Polish space with a~$\zero'$-computable compact presentation, then it has a low Polish presentation. We show that it fails for some non-perfect Polish space.
\begin{proposition}\label{prop_nonperfect}
There exists a (non-perfect) compact Polish space~$\X$ with a~$\zero'$-computable compact presentation, but no low Polish presentation.
\end{proposition}

\begin{proof}
By Theorem \ref{thm:Cantor-Bendixson-CP0} applied to~$\d=\zero'$, there is a space~$\X$ having a~$\zero'$-computable compact presentation, but such that its Cantor-Bendixson derivative~$\X'$ does not have a~$\zero''$-computable compact presentation. If~$\X$ had a low Polish presentation then by Lemma \ref{lem:CB-derivative},~$\X'$ would have a~$\zero''$-computable compact presentation.
\end{proof}

Therefore, for~$\d=\zero'$, this space~$\X$ has a~$\d$-computable compact presentation, but there is no~$\e$ such that~$X$ has an~$\e$-computable Polish presentation and~$\e'\leq \d\oplus \zero'=\zero'$.


\end{document}